\documentclass[reqno,11pt]{amsart}
\usepackage{amssymb}
\usepackage{graphicx}
\usepackage{pstricks}
\usepackage{amsmath}
\usepackage{amsmath}
\usepackage{amsxtra}

\theoremstyle{plain}
\newtheorem{theorem}{Theorem}[section]
\newtheorem{lemma}[theorem]{Lemma}

\newtheorem{proposition}[theorem]{Proposition}
\newtheorem{corollary}[theorem]{Corollary}

\newtheorem{problem}[theorem]{Problem}

\theoremstyle{definition}

\newtheorem{example}[theorem]{Example}
\newtheorem{remark}[theorem]{Remark}

\newtheorem*{acknowledgements}{Acknowledgement}

\numberwithin{equation}{section} 
\begin{document}
\title{When is hyponormality for $2$-variable weighted shifts invariant
under powers?}
\author{Ra\'{u}l E. Curto}
\address{Department of Mathematics, The University of Iowa, Iowa City, Iowa
52242}
\email{raul-curto@uiowa.edu}
\urladdr{http://www.math.uiowa.edu/\symbol{126}rcurto/}
\author{Jasang Yoon}
\address{Department of Mathematics, The University of Texas-Pan American,
Edinburg, Texas 78539}
\email{yoonj@utpa.edu}
\urladdr{http://www.math.utpa.edu/\symbol{126}yoonj/}
\thanks{The first named author was partially supported by NSF Grants
DMS-0400741 and DMS-0801168.}
\thanks{The second named author was partially supported by a Faculty
Research Council Grant at The University of Texas-Pan American.}

\begin{abstract}
For $2$-variable weighted shifts $W_{(\alpha ,\beta )}\equiv (T_{1},T_{2})$
we study the invariance of (joint) $k$-hyponormality under the action $%
(h,\ell )\mapsto W_{(\alpha ,\beta )}^{(h,\ell )}:=(T_{1}^{h},T_{2}^{\ell
})\;\;(h,\ell \geq 1)$. \ We show that for every $k\geq 1$ there exists $%
W_{(\alpha ,\beta )}$ such that $W_{(\alpha ,\beta )}^{(h,\ell )}$ is $k$%
-hyponormal (all $h\geq 2$, $\ell \geq 1$) but $W_{(\alpha ,\beta )}$ is not
$k$-hyponormal. \ On the positive side, for a class of $2$-variable weighted
shifts with tensor core we find a computable necessary condition for
invariance. \ Next, we exhibit a large nontrivial class for which
hyponormality is indeed invariant under \textit{all} powers; moreover, for
this class $2$-hyponormality automatically implies subnormality. \ Our results partially
depend on new formulas for the determinant of generalized Hilbert matrices
and on criteria for their positive semi-definiteness.
\end{abstract}

\maketitle

\section{\label{Int}Introduction}

Given a pair $\mathbf{T}\equiv (T_{1},T_{2})$ of commuting subnormal Hilbert
space operators, the Lifting Problem for Commuting Subnormals (LPCS) calls
for necessary and sufficient conditions for the existence of a commuting
pair $\mathbf{N}\equiv (N_{1},N_{2})$ of normal extensions of $T_{1}$ and $%
T_{2}$. \ In previous work (\cite{CLY1}, \cite{CLY2}, \cite{ROMP}, \cite%
{CLY4}, \cite{CuYo1}, \cite{CuYo2}, \cite{CuYo3}) we have studied the
relevance of (joint) $k$-hyponormality to LPCS. \ In particular, one asks to
what extent the existence of liftings for the powers $\mathbf{T}^{(h,\ell
)}\equiv (T_{1}^{h},T_{2}^{\ell })\;\;(h,\ell \geq 1)$ can guarantee a
lifting for $\mathbf{T}$. \ For the class of $2$-variable weighted shifts $%
W_{(\alpha ,\beta )}$, it is often the case that the powers are less complex
than the initial pair; thus it becomes especially significant to unravel the
invariance of $k$-hyponormality under the action $(h,\ell )\mapsto
W_{(\alpha ,\beta )}^{(h,\ell )}\;\;(h,\ell \geq 1)$. \

Our aim in this paper is to shed new light on some of the intricacies
associated with LPCS and $k$-hyponormality for powers of commuting
subnormals. \ To describe our results we need some notation; we further
expand on our terminology and basic results in Section \ref{Sect2}. \ We use
$\mathfrak{H}_{0}$ (resp. $\mathfrak{H}_{\infty })$ to denote the set of
commuting pairs of subnormal operators (resp. subnormal pairs) on Hilbert
space. \ For $k\geq 1$, we let $\mathfrak{H}_{k}$ denote the class of $k$%
-hyponormal pairs in $\mathfrak{H}_{0}$. \ Clearly, $\mathfrak{H}_{\infty
}\subseteq \cdots \subseteq \mathfrak{H}_{k}\subseteq \cdots \subseteq
\mathfrak{H}_{2}\subseteq \mathfrak{H}_{1}\subseteq \mathfrak{H}_{0}$. \ The
main results in \cite{CuYo1} and \cite{CLY1} show that these inclusions are
all proper. \ In our previous research we have shown that detecting these
proper inclusions can be done within classes of $2$-variable weighted shifts
with relatively simple weight structure, as we now describe. \

For a sequence $\alpha \equiv \{\alpha _{k}\}_{k=0}^{\infty }\in \ell
^{\infty }(\mathbb{Z}_{+})$ of positive numbers, we let $W_{\alpha }\equiv
\operatorname{shift}\;(\alpha _{0},\alpha _{1},\cdots )$ denote the unilateral weighted shift
on $\ell ^{2}(\mathbb{Z}_{+})$ given by $W_{\alpha }e_{k}:=\alpha
_{k}e_{k+1}\;(k\geq 0)$. \ We also let $U_{+}:=\operatorname{shift}\;(1,1,\cdots )$
(the (unweighted) unilateral shift), and for $0<a<1$ we let $S_{a}:=\operatorname{%
shift}\;(a,1,1,\cdots )$. \ Multivariable weighted shifts are defined in an
analogous manner. \ For instance, on $\ell ^{2}(\mathbb{Z}_{+}^{2})$ we let $%
W_{(\alpha ,\beta )}\equiv (T_{1},T_{2})$ denote the $2$-variable weighted
shift associated with weight sequences $\alpha $ and $\beta $, defined by $%
T_{1}e_{\mathbf{k}}:=\alpha _{\mathbf{k}}e_{\mathbf{k}+\mathbf{\varepsilon }%
_{1}}$ and $T_{2}e_{\mathbf{k}}:=\beta _{\mathbf{k}}e_{\mathbf{k}+\mathbf{%
\varepsilon }_{2}}\;(\mathbf{k}\in \mathbb{Z}_{+}^{2})$. \

For an arbitrary $2$-variable weighted shift $W_{(\alpha ,\beta )}$, we let $%
\mathcal{M}_{i}$ (resp. $\mathcal{N}_{j})$ be the subspace of $\ell ^{2}(%
\mathbb{Z}_{+}^{2})$ which is spanned by the canonical orthonormal basis
associated to indices $\mathbf{k}=(k_{1},k_{2})$ with $k_{1}\geq 0$ and $%
k_{2}\geq i$ $($resp. $k_{1}\geq j$ and $k_{2}\geq 0$). \ We will often
write $\mathcal{M}_{1}$ simply as $\mathcal{M}$ and $\mathcal{N}_{1}$ as $%
\mathcal{N}$. \ The core $c(W_{(\alpha ,\beta )})$ of $W_{(\alpha ,\beta )}$
is the restriction of $W_{(\alpha ,\beta )}$ to the invariant subspace $%
\mathcal{M}\bigcap \mathcal{N}$. \ A $2$-variable weighted shift $W_{(\alpha
,\beta )}$ is said to be \textit{of tensor form} if it is of the form $%
(I\otimes W_{\alpha },W_{\beta }\otimes I)$. \ The class of all $2$-variable
weighted shifts $W_{(\alpha ,\beta )}\in \mathfrak{H}_{0}$ whose core is of
tensor form will be denoted by $\mathcal{TC}$; in symbols, $\mathcal{TC}%
:=\{W_{(\alpha ,\beta )}\in \mathfrak{H}_{0}:c(W_{(\alpha ,\beta )})\text{
is of tensor form}\}$.\ \

We now consider the class $\mathcal{S}:=\{W_{(\alpha ,\beta )}\in \mathfrak{H%
}_{0}:\alpha _{(k_{1},0)}=\alpha _{(k_{1}+1,0)}$ and $\beta
_{(0,k_{2})}=\beta _{(0,k_{2}+1)}$ for some $k_{1}\geq 1$ and $k_{2}\geq 1\}$
and we let $\mathcal{S}_{1}:=\mathcal{S}\cap \mathfrak{H}_{1}$. \ From
propagation phenomena for $1$- and $2$-variable weighted shifts (see \cite%
{CuYo2}, \cite{CLY4}), we observe that, without loss of generality, we
can always assume that the restrictions of each $W_{(\alpha ,\beta )}\in
\mathcal{S}_{1}$ to the invariant subspace $\mathcal{M}$ (resp. $\mathcal{N}$%
) is of the form $\left( I\otimes S_{a},U_{+}\otimes I\right) $ (resp. $%
\left( I\otimes U_{+},S_{b}\otimes I\right) $); cf. Figure \ref{S1}(i). \
In particular, the core $c(W_{(\alpha ,\beta )})$ of a $2$-variable weighted
shift in $\mathcal{S}_{1}$ is always the doubly commuting pair $\left(
I\otimes U_{+},U_{+}\otimes I\right) $; as a result, $W_{(\alpha ,\beta
)}\in \mathcal{TC}$. \ Observe also that if $W_{(\alpha ,\beta )}\in
\mathcal{S}_{1}$, then $W_{(\alpha ,\beta )}$ is completely determined by
the three parameters $x:=\alpha _{(0,0)}$, $y:=\beta _{(0,0)}$ and $%
a:=\alpha _{(0,1)}$. $\ $Thus we shall often denote a $2$-variable weighted
shift $W_{(\alpha ,\beta )}\in \mathcal{S}_{1}$ by $\left\langle
x,y,a\right\rangle $. \

Between $\mathcal{S}_{1}$ and $\mathcal{TC}$ there is a class that provides
significant information about LPCS, and we now define it. \ Let $\mathcal{A}%
:=\{W_{(\alpha ,\beta )}\in \mathcal{TC}:c(W_{(\alpha ,\beta )})$ is $1$%
-atomic$\}$. \ Clearly $\mathcal{S}_{1}\varsubsetneq \mathcal{A}%
\varsubsetneq \mathcal{TC}$. \ In \cite{ROMP} we solved LPCS within the
class $\mathcal{TC}$, and in particular we gave a simple test for
subnormality within $\mathcal{A}$. \

To prove that the $k$-hyponormality of all powers need not guarantee the $k$%
-hyponormality of the initial pair, we build an example that uses weights
related to those of the Bergman shift. \ The reader will recall that the
moment matrix associated with the Bergman shift is the classical Hilbert
matrix. \ Thus to deal with our situation we need to describe positivity and
the calculation of determinants for \textit{generalized} Hilbert matrices;
we do this in Theorem \ref{Thm-Hilbert}. \ Although Section \ref{GHM} has
intrinsic and independent value since it deals with matrices that arise in
various contexts, the main reason for including it here is that it contributes a basic tool for
producing some of the examples in subsequent sections.

It is well known that for a general operator $T$ on Hilbert space, the
hyponormality of $T$ does not imply the hyponormality of $T^{2}$ \cite{Hal}.
\ However, for a unilateral weighted shift $W_{\alpha }$, the hyponormality
of $W_{\alpha }$ (detected by the condition $\alpha _{j}\leq \alpha _{j+1}$
for all $j\geq 0$) does imply the hyponormality of every power $W_{\alpha
}^{n}\;(n\geq 2)$. \ It is also well known that the subnormality of $T$
implies the subnormality of $T^{n}$ (all $n\geq 2$), but the converse
implication is not true, even if $T$ is a unilateral weighted shift \cite%
{Sta2}. \ Since $k$-hyponormality lies between hyponormality and
subnormality, it is then natural to consider

\begin{problem}
\label{question0}Let $T$ be an operator and let $k\geq 2$. \ \newline
(i) \ Does the $k$-hyponormality of $T$ imply the $k$-hyponormality of $%
T^{2} $?\newline
(ii) \ Does the $k$-hyponormality of $T^{2}$ imply the $k$-hyponormality of $%
T$?
\end{problem}

At the beginning of Section \ref{Sect4} we consider this problem, and we
subsequently study its multivariable analogue. \ It is worth noting that, in the
multivariable case, the standard assumption on a pair $\mathbf{T}\equiv
(T_{1},T_{2})$ is that each component $T_{i}$ be subnormal ($i=1,2$). \ With
this in mind, comparing the $k$-hyponormality of a $2$-variable weighted
shift $W_{(\alpha ,\beta )}\in \mathfrak{H}_{0}$ to the $k$-hyponormality of
its powers $W_{(\alpha ,\beta )}^{(h,\ell )}$ is highly nontrivial. \ We now
formulate the relevant problems in the multivariable case.

\begin{problem}
\label{question2}Given $k\geq 1$ and $W_{(\alpha ,\beta )}\in \mathfrak{H}%
_{k}$, does it follow that $W_{(\alpha ,\beta )}^{(h,\ell )}\in \mathfrak{H}%
_{k}$ for all $h,\ell \geq 1$?
\end{problem}

In Section \ref{hypoinv} we establish that subclasses of the class $\mathfrak{H}_{k}$ $%
(k\geq 1)$ are often invariant under powers. \ Concretely, we prove
that there exists a rich collection of $2$-variable weighted shifts $W_{(\alpha ,\beta )}\in \mathfrak{H}_{2}$ such that $%
W_{(\alpha ,\beta )}^{(2,1)}\in \mathfrak{H}_{2}$ (Theorem \ref{Thm8}). \
Conversely, we can ask

\begin{problem}
\label{question1}Given $k\geq 1$, assume that for all $h\geq 2$ and $\ell
\geq 1$, $W_{(\alpha ,\beta )}^{(h,\ell )}\in \mathfrak{H}_{k}$. \ Does it
follow that $W_{(\alpha ,\beta )}\in \mathfrak{H}_{k}$?
\end{problem}

In Theorem \ref{Thm(invariant)-00} we answer Problem \ref{question1} in the
negative; that is, for each $k\geq 1$ we build a $2$-variable weighted shift
$W_{(\alpha ,\beta )}\in \mathfrak{H}_{0}\;\backslash \;\mathfrak{H}_{k}$
such that $W_{(\alpha ,\beta )}^{(h,\ell )}\in \mathfrak{H}_{k}$ (all $h\geq
2$ and $\ell \geq 1$). \

Next, for $k=1,2$, we find a computable necessary condition for the $k$%
-hyponormality of $W_{(\alpha ,\beta )}$ to remain invariant under all
powers (Theorem \ref{Thm-main}). \ We then show that this necessary
conditions is not sufficient (Remark \ref{Re4}(ii)). \

Section \ref{S1inv} is devoted to the study of the class $\mathcal{S}_{1}$.
\ We show that for $\left\langle x,y,a\right\rangle \in \mathcal{S}_{1}$,
all powers $\left\langle x,y,a\right\rangle ^{(h,\ell )}$ are hyponormal
(Theorem \ref{hyponormal}). \ Moreover, a shift $\left\langle
x,y,a\right\rangle \in \mathcal{S}_{1}$ is $2$-hyponormal if and only if it
is subnormal.

As we mentioned before, for single operators it is an open problem whether
the $2$-hyponormality of $T$ implies the $2$-hyponormality of $T^{2}$. \
Although this problem is intimately related to Theorem \ref{Thm8}, we
observe that the latter does not provide an answer to Problem \ref{question0}
when $k=2$, since our pairs consist of commuting subnormal operators.

Problem \ref{question2} is a special case of a much more general problem,
that of determining necessary and sufficient conditions for the weak $k$%
-hyponormality of a commuting pair. \ We say that a pair $\mathbf{T}\in
\mathfrak{H}_{0}$ is \textit{weakly }$k$\textit{-hyponormal} if
\begin{equation*}
\mathbf{p}(\mathbf{T}):=(p_{1}(T_{1},T_{2}),p_{2}(T_{1}T_{2})),
\end{equation*}%
is hyponormal for all polynomials $p_{1}$, $p_{2}\in \mathbb{C}[z,w]$ with $%
\deg \;p_{1},\deg \;p_{2}\leq k$, where $\mathbf{p}\equiv (p_{1},p_{2})$. \
To verify that $\mathbf{T}$ is weakly $k$-hyponormality is highly
nontrivial. \ Thus Problems \ref{question2} and \ref{question1} can be
regarded as suitably multivariable analogues of \cite[Question 33]{Shi}: If $%
T$ is a hyponormal unilateral shift and if $p$ is a polynomial, must $p(T)$
be hyponormal? \ If $T$ is subnormal, the answer is clearly yes, but we note
that polynomial hyponormality is strictly weaker than subnormality, as
proved in \cite{CuPu}.

\begin{acknowledgements}
Most of the examples and several proofs in this paper were obtained using
calculations with the software tool \textit{Mathematica \cite{Wol}.}
\end{acknowledgements}

%%%%%%%%%%%%%%%%%%%%%%

\section{\label{Sect2}Notation and Preliminaries}

Let $\mathcal{H}$ be a complex Hilbert space and let $\mathcal{B}(\mathcal{H}%
)$ denote the algebra of bounded linear operators on $\mathcal{H}$. $\ $For $%
S,T\in \mathcal{B}(\mathcal{H})$ let $[S,T]:=ST-TS$. \ We say that an $n$%
-tuple $\mathbf{T}=(T_{1},\cdots ,T_{n})$ of operators on $\mathcal{H}$ is
(jointly) \textit{hyponormal} if the operator matrix
\begin{equation*}
\lbrack \mathbf{T}^{\ast },\mathbf{T]:=}\left(
\begin{array}{llll}
\lbrack T_{1}^{\ast },T_{1}] & [T_{2}^{\ast },T_{1}] & \cdots & [T_{n}^{\ast
},T_{1}] \\
\lbrack T_{1}^{\ast },T_{2}] & [T_{2}^{\ast },T_{2}] & \cdots & [T_{n}^{\ast
},T_{2}] \\
\text{ \thinspace \thinspace \quad }\vdots & \text{ \thinspace \thinspace
\quad }\vdots & \ddots & \text{ \thinspace \thinspace \quad }\vdots \\
\lbrack T_{1}^{\ast },T_{n}] & [T_{2}^{\ast },T_{n}] & \cdots & [T_{n}^{\ast
},T_{n}]%
\end{array}%
\right)
\end{equation*}%
is positive on the direct sum of $n$ copies of $\mathcal{H}$ (cf. \cite{Ath}%
, \cite{CMX}). \ The $n$-tuple $\mathbf{T}$ is said to be \textit{normal} if
$\mathbf{T}$ is commuting and each $T_{i}$ is normal, and $\mathbf{T}$ is
\textit{subnormal }if $\mathbf{T}$ is the restriction of a normal $n$-tuple
to a common invariant subspace. \ For $k\geq 1$, a commuting pair $\mathbf{T}%
\equiv (T_{1},T_{2})$ is said to be $k$\textit{-hyponormal} (\cite{CLY1}) if
\begin{equation*}
\mathbf{T}(k):=(T_{1},T_{2},T_{1}^{2},T_{2}T_{1},T_{2}^{2},\cdots
,T_{1}^{k},T_{2}T_{1}^{k-1},\cdots ,T_{2}^{k})
\end{equation*}%
is hyponormal, or equivalently
\begin{equation*}
\lbrack \mathbf{T}(k)^{\ast },\mathbf{T}(k)]=([(T_{2}^{q}T_{1}^{p})^{\ast
},T_{2}^{m}T_{1}^{n}])_{_{1\leq p+q\leq k}^{1\leq n+m\leq k}}\geq 0.
\end{equation*}%
Clearly, normal $\Rightarrow $ subnormal $\Rightarrow $ $k$-hyponormal. \
The Bram-Halmos criterion states that an operator $T\in \mathcal{B}(\mathcal{%
H})$ is subnormal if and only if the $k$-tuple $(T,T^{2},\cdots ,T^{k})$ is
hyponormal for all $k\geq 1$.

For $\alpha \equiv \{\alpha _{n}\}_{n=0}^{\infty }$ a bounded sequence of
positive real numbers (called \textit{weights}), let $W_{\alpha }:\ell ^{2}(%
\mathbb{Z}_{+})\rightarrow \ell ^{2}(\mathbb{Z}_{+})$ be the associated
unilateral weighted shift, defined by $W_{\alpha }e_{n}:=\alpha
_{n}e_{n+1}\; $(all $n\geq 0$), where $\{e_{n}\}_{n=0}^{\infty }$ is the
canonical orthonormal basis in $\ell ^{2}(\mathbb{Z}_{+})$. \ The moments of
$\alpha $ are given as
\begin{equation*}
\gamma _{k}\equiv \gamma _{k}(\alpha ):=\left\{
\begin{array}{cc}
1, & \text{if }k=0 \\
\alpha _{0}^{2}\cdots \alpha _{k-1}^{2}, & \text{if }k>0%
\end{array}%
\right\} .
\end{equation*}%
It is easy to see that $W_{\alpha }$ is never normal, and that it is
hyponormal if and only if $\alpha _{0}\leq \alpha _{1}\leq \cdots $. \
Similarly, consider double-indexed positive bounded sequences $\alpha _{%
\mathbf{k}},\beta _{\mathbf{k}}\in \ell ^{\infty }(\mathbb{Z}_{+}^{2})$, $%
\mathbf{k}\equiv (k_{1},k_{2})\in \mathbb{Z}_{+}^{2}:=\mathbb{Z}_{+}\times
\mathbb{Z}_{+}$ and let $\ell ^{2}(\mathbb{Z}_{+}^{2})$\ be the Hilbert
space of square-summable complex sequences indexed by $\mathbb{Z}_{+}^{2}$.
\ (Recall that $\ell ^{2}(\mathbb{Z}_{+}^{2})$ is canonically isometrically
isomorphic to $\ell ^{2}(\mathbb{Z}_{+})\bigotimes \ell ^{2}(\mathbb{Z}_{+})$%
.) \ We define the $2$-variable weighted shift $W_{(\alpha ,\beta )}\equiv
(T_{1},T_{2})$\ by
\begin{equation*}
T_{1}e_{\mathbf{k}}:=\alpha _{\mathbf{k}}e_{\mathbf{k+}\varepsilon _{1}}
\end{equation*}%
\begin{equation*}
T_{2}e_{\mathbf{k}}:=\beta _{\mathbf{k}}e_{\mathbf{k+}\varepsilon _{2}},
\end{equation*}%
where $\mathbf{\varepsilon }_{1}:=(1,0)$ and $\mathbf{\varepsilon }%
_{2}:=(0,1)$. \ Clearly,
\begin{equation}
T_{1}T_{2}=T_{2}T_{1}\Longleftrightarrow \beta _{\mathbf{k+}\varepsilon
_{1}}\alpha _{\mathbf{k}}=\alpha _{\mathbf{k+}\varepsilon _{2}}\beta _{%
\mathbf{k}}\;(\text{all }\mathbf{k}\in \mathbb{Z}_{+}^{2}).
\label{commuting}
\end{equation}%
In an entirely similar way one can define multivariable weighted shifts. \
Trivially, a pair of unilateral weighted shifts $W_{a}$ and $W_{\beta }$
gives rise to a $2$-variable weighted shift $W_{(\alpha ,\beta )}\equiv
(T_{1},T_{2})$, if we let $\alpha _{(k_{1},k_{2})}:=\alpha _{k_{1}}$ and $%
\beta _{(k_{1},k_{2})}:=\beta _{k_{2}}\;($all $k_{1},k_{2}\in \mathbb{Z}%
_{+}) $. \ In this case, $W_{(\alpha ,\beta )}$ is subnormal (resp.
hyponormal) if and only if so are $T_{1}$ and $T_{2}$; in fact, under the
canonical identification of $\ell ^{2}(\mathbb{Z}_{+}^{2})$ with $\ell ^{2}(%
\mathbb{Z}_{+})\bigotimes \ell ^{2}(\mathbb{Z}_{+})$, we have $T_{1}\cong
I\bigotimes W_{a}$ and $T_{2}\cong W_{\beta }\bigotimes I$, and $W_{(\alpha
,\beta )}$ is also doubly commuting. \ For this reason, we do not focus
attention on shifts of this type, and use them only when the above mentioned
triviality is desirable or needed.

Given $\mathbf{k}\in \mathbb{Z}_{+}^{2}$, the moment of $(\alpha ,\beta )$
of order $\mathbf{k}$ is
\begin{equation*}
\gamma _{\mathbf{k}}\equiv \gamma _{\mathbf{k}}(\alpha ,\beta ):=%
\begin{cases}
1, & \text{if }\mathbf{k}=\mathbf{0} \\
\alpha _{(0,0)}^{2}\cdots \alpha _{(k_{1}-1,0)}^{2}, & \text{if }k_{1}\geq 1%
\text{ and }k_{2}=0 \\
\beta _{(0,0)}^{2}\cdots \beta _{(0,k_{2}-1)}^{2}, & \text{if }k_{1}=0\text{
and }k_{2}\geq 1 \\
\alpha _{(0,0)}^{2}\cdots \alpha _{(k_{1}-1,0)}^{2}\beta
_{(k_{1},0)}^{2}\cdots \beta _{(k_{1},k_{2}-1)}^{2}, & \text{if }k_{1}\geq 1%
\text{ and }k_{2}\geq 1.%
\end{cases}%
\end{equation*}%
We remark that, due to the commutativity condition (\ref{commuting}), $%
\gamma _{\mathbf{k}}$ can be computed using any nondecreasing path from $%
(0,0)$ to $(k_{1},k_{2})$. \

We now recall a well known characterization of subnormality for
multivariable weighted shifts \cite{JeLu}, which in the single variable case is due to C. Berger (cf. \cite[%
III.8.16]{Con}) and was independently established by R. Gellar and L.J. Wallen
\cite{GeWa}: $\ W_{(\alpha ,\beta )}$ admits a
commuting normal extension if and only if there is a probability measure $%
\mu $ (which we call the \textit{Berger measure} of $W_{(\alpha ,\beta )}$)
defined on the $2$-dimensional rectangle $R=[0,a_{1}]\times \lbrack 0,a_{2}]$
(where $a_{i}:=\left\| T_{i}\right\| ^{2}$) such that $\gamma _{\mathbf{k}%
}=\int_{R}\mathbf{t}^{\mathbf{k}}d\mu (\mathbf{t}):=%
\int_{R}t_{1}^{k_{1}}t_{2}^{k_{2}}d\mu (\mathbf{t}),$ for all $\mathbf{k}\in
\mathbb{Z}_{+}^{2}$. \ Observe that $U_{+}$ and $S_{a}$ are subnormal, with
Berger measures $\delta _{1}$ and $(1-a^{2})\delta _{0}+a^{2}\delta _{1} $,
respectively, where $\delta _{p}$ denotes the point-mass probability measure
with support the singleton set $\{p\}$. \ Also, a $2$-variable weighted
shift $W_{(\alpha ,\beta )}\in \mathcal{S}_{1}$ has a core with Berger
measure $\delta _{1}\times \delta _{1}$.

%%%%%%%%%%%%%%%%%%%%%%%%

\section{\label{GHM}The determinant of a generalized Hilbert matrix}

Given positive real numbers $x$ and $h$, and an integer $k\geq 1$, we define
the generalized Hilbert matrix $A_{k}(x,h)$ as follows:%
\begin{equation*}
(A_{k}(x,h))_{i,j}:=\left\{
\begin{tabular}{ll}
$x,$ & if $i=j=1$ \\
$\frac{1}{(i+j-2)h+1},$ & otherwise%
\end{tabular}%
\right. \text{ \ \ \ \ }(1\leq i,j\leq k+1).
\end{equation*}%
(Observe that $A_{k}(1,1)$ is the classical Hilbert matrix.) \ In this
section we calculate the determinant of, and establish positivity properties
for, the generalized Hilbert matrix $A_{k}(x,h)$.\

To describe our results, we need some notation. \ We let $0!:=1$, $%
k!:=k(k-1)!$, and $k^{!}:=\Pi _{i=1}^{k}i!$. \ We also let
\begin{eqnarray*}
f_{0} &:&=x, \\
f_{\ell +1} &:&=f_{\ell }\left( \frac{(k-\ell )h+1}{\left( k-\ell \right) h}%
\right) ^{2}-\frac{2(k-\ell )h+1}{\left( \left( k-\ell \right) h\right) ^{2}}%
\;(0\leq \ell \leq k-1),
\end{eqnarray*}%
\begin{equation}
\begin{tabular}{l}
$g(h,k):=\left( \frac{1}{(k+1)h+1}\right) ^{k}\prod_{j=0}^{k-1}\left( \frac{1%
}{(jh+1)(2kh-jh+1)}\right) ^{j+1},$%
\end{tabular}
\label{gk}
\end{equation}%
and%
\begin{equation}
\begin{tabular}{l}
$f(x,h,k):=x\left( \frac{kh+1}{kh}\right) ^{2}\left( \frac{\left( k-1\right)
h+1}{\left( k-1\right) h}\right) ^{2}\cdots \left( \frac{3h+1}{3h}\right)
^{2}\left( \frac{2h+1}{2h}\right) ^{2}\left( \frac{(h+1)^{2}}{h}\right) $ \\
$\ \ \ \ \ \ \ \ \ \ \ \ \ \ -\left( \frac{2\left( kh\right) +1}{\left(
kh\right) ^{2}}\right) \left( \frac{\left( k-1\right) h+1}{\left( k-1\right)
h}\right) ^{2}\cdots \left( \frac{3h+1}{3h}\right) ^{2}\left( \frac{2h+1}{2h}%
\right) ^{2}\frac{(h+1)^{2}}{h}$ \\
$\ \ \ \ \ \ \ \ \ \ \ \ \ \ -\left( \frac{2\left( \left( k-1\right)
h\right) +1}{\left( \left( k-1\right) h\right) ^{2}}\right) \left( \frac{%
\left( k-2\right) h+1}{\left( k-2\right) h}\right) ^{2}\cdots \left( \frac{%
3h+1}{3h}\right) ^{2}\left( \frac{2h+1}{2h}\right) ^{2}\frac{(h+1)^{2}}{h}$
\\
$\ \ \ \ \ \ \ \ \ \ \ \ \ \ -\cdots -\left( \frac{2\left( 2h\right) +1}{%
\left( 2h\right) ^{2}}\right) \frac{(h+1)^{2}}{h}-\left( \frac{2h+1}{h}%
\right) .$%
\end{tabular}
\label{fx}
\end{equation}

\begin{theorem}
\label{Thm-Hilbert}For $x,h>0$ and $k\geq 1$, we have
\begin{equation}
\det \;A_{k}(x,h)=h^{k(k+1)}\left( k^{!}\right) ^{2}g(h,k)f_{k},
\label{eqHilbert}
\end{equation}%
where $f_{k}=f_{k-1}\left( \frac{h+1}{h}\right) ^{2}-\frac{2h+1}{h^{2}}$. \
Moreover,
\begin{equation*}
f_{k}=f(x,h,k).
\end{equation*}
\end{theorem}

\begin{proof}
Consider the $(k+1)\times (k+1)$ matrix%
\begin{equation*}
A_{k}(x,h)=\left(
\begin{array}{cccccc}
x & \frac{1}{h+1} & \frac{1}{2h+1} & \cdots & \frac{1}{(k-1)h+1} & \frac{1}{%
kh+1} \\
\frac{1}{h+1} & \frac{1}{2h+1} & \frac{1}{3h+1} & \cdots & \frac{1}{kh+1} &
\frac{1}{(k+1)h+1} \\
\frac{1}{2h+1} & \frac{1}{3h+1} & \frac{1}{4h+1} & \cdots & \frac{1}{(k+1)h+1%
} & \frac{1}{(k+2)h+1} \\
\vdots & \vdots & \vdots & \ddots & \vdots & \vdots \\
\frac{1}{(k-1)h+1} & \frac{1}{kh+1} & \frac{1}{(k+1)h+1} & \cdots & \frac{1}{%
(2k-2)h+1} & \frac{1}{(2k-1)h+1} \\
\frac{1}{kh+1} & \frac{1}{(k+1)h+1} & \frac{1}{(k+2)h+1} & \cdots & \frac{1}{%
(2k-1)h+1} & \frac{1}{2kh+1}%
\end{array}%
\right) .
\end{equation*}%
Let us first subtract the $(k+1)$-st row from each row above it. \ The entry
in the $j$-th column of the $i$-th row becomes%
\begin{equation*}
\left\{
\begin{tabular}{ll}
$x-\frac{1}{kh+1},$ & if $(i,j)=(1,1)$ \\
$\frac{1}{(i+j-2)h+1}-\frac{1}{(k+j-1)h+1}=\frac{\left( k-i+1\right) h}{%
\left[ (i+j-2)h+1\right] \left[ (k+j-1)h+1\right] },$ & if $(i,j)\neq (1,1).$%
\end{tabular}%
\right.
\end{equation*}%
The new $(k+1)\times (k+1)$ matrix is%
\begin{equation*}
B_{k}(x,h):=\left(
\begin{array}{ccccc}
x-\frac{1}{kh+1} & \frac{kh}{\left[ h+1\right] \left[ (k+1)h+1\right] } &
\cdots & \frac{kh}{\left[ (k-1)h+1\right] \left[ (2k-1)h+1\right] } & \frac{%
kh}{\left[ kh+1\right] \left[ 2kh+1\right] } \\
\frac{\left( k-1\right) h}{\left[ h+1\right] \left[ kh+1\right] } & \frac{%
\left( k-1\right) h}{\left[ 2h+1\right] \left[ (k+1)h+1\right] } & \cdots &
\frac{\left( k-1\right) h}{\left[ kh+1\right] \left[ (2k-1)h+1\right] } &
\frac{\left( k-1\right) h}{\left[ (k+1)h+1\right] \left[ 2kh+1\right] } \\
\vdots & \vdots & \ddots & \vdots & \vdots \\
\frac{h}{\left[ (k-1)h+1\right] \left[ kh+1\right] } & \frac{h}{\left[ kh+1%
\right] \left[ (k+1)h+1\right] } & \cdots & \frac{h}{\left[ (2k-2)h+1\right] %
\left[ (2k-1)h+1\right] } & \frac{h}{\left[ (2k-1)h+1\right] \left[ 2kh+1%
\right] } \\
\frac{1}{kh+1} & \frac{1}{(k+1)h+1} & \cdots & \frac{1}{(2k-1)h+1} & \frac{1%
}{2kh+1}%
\end{array}%
\right) .
\end{equation*}%
Note that $\det \;A_{k}(x,h)=\det \;B_{k}(x,h)$. \ To compute $\det
\;B_{k}(x,h)$, we observe that one can factor out $(k-(i-1))h$ from the $i$%
-th row ($1\leq i<k$) and $\frac{1}{(k+j-1)h+1}$ from the $j$-th column ($%
1\leq j\leq k+1$) in the matrix $B_{k}(x,h)$. \ Hence we obtain
\begin{equation*}
\det \;A_{k}(x,h)=k!h^{k}\cdot \frac{1}{kh+1}\cdot \frac{1}{(k+1)h+1}\cdots
\frac{1}{2kh+1}\cdot \det \;C_{k}(x,h),
\end{equation*}%
where%
\begin{equation*}
C_{k}(x,h):=\left(
\begin{array}{cccc}
\left( x-\frac{1}{kh+1}\right) \frac{kh+1}{kh} & \frac{1}{h+1} & \cdots &
\frac{1}{kh+1} \\
\frac{1}{h+1} & \frac{1}{2h+1} & \cdots & \frac{1}{(k+1)h+1} \\
\vdots & \vdots & \ddots & \vdots \\
\frac{1}{(k-1)h+1} & \frac{1}{kh+1} & \cdots & \frac{1}{(2k-1)h+1} \\
1 & 1 & \cdots & 1%
\end{array}%
\right) .
\end{equation*}%
Next, let us subtract the last column from each of the preceding columns in
the $(k+1)\times (k+1)$ matrix $C_{k}(x,h)$. \ We obtain%
\begin{equation*}
D_{k}(x,h):=\left(
\begin{array}{ccccc}
\left( x-\frac{1}{kh+1}\right) \frac{kh+1}{kh}-\frac{1}{kh+1} & \frac{\left(
k-1\right) h}{(h+1)\left[ kh+1\right] } & \cdots & \frac{h}{\left[ (k-1)h+1%
\right] \left[ kh+1\right] } & \frac{1}{kh+1} \\
\frac{kh}{(h+1)\left[ (k+1)h+1\right] } & \frac{\left( k-1\right) h}{(2h+1)%
\left[ (k+1)h+1\right] } & \cdots & \frac{h}{\left[ kh+1\right] \left[
(k+1)h+1\right] } & \frac{1}{(k+1)h+1} \\
\vdots & \vdots & \ddots & \vdots & \vdots \\
\frac{kh}{\left[ (k-1)h+1\right] \left[ (2k-1)h+1\right] } & \frac{\left(
k-1\right) h}{\left[ kh+1\right] \left[ (2k-1)h+1\right] } & \cdots & \frac{h%
}{\left[ (2k-2)h+1\right] \left[ (2k-1)h+1\right] } & \frac{1}{(2k-1)h+1} \\
0 & 0 & \cdots & 0 & 1%
\end{array}%
\right) .
\end{equation*}%
Note that $\det \;C_{k}(x,h)=\det \;D_{k}(x,h)$. \ As we have done before,
let us factor out $(k-(j-1))h$ from the $j$-th column ($1\leq j\leq k$) and $%
\frac{1}{(k+i-1)h+1}$ from the $i$-th row ($1\leq i\leq k$) in the matrix $%
D_{k}(x,h)$. \ Let
\begin{equation*}
f_{1}\equiv f_{1}(x,h,k):=x\left( \frac{kh+1}{kh}\right) ^{2}-\frac{2kh+1}{%
\left( kh\right) ^{2}}.
\end{equation*}%
Then we have
\begin{equation*}
\begin{tabular}{l}
$\det \;A_{k}(x,h)=\left( k!\right) ^{2}h^{2k}\left( \frac{1}{kh+1}\right)
^{2}\left( \frac{1}{(k+1)h+1}\right) ^{2}\cdots \left( \frac{1}{(2k-1)h+1}%
\right) ^{2}\left( \frac{1}{2kh+1}\right) \det \;A_{k-1}(f_{1},h),$%
\end{tabular}%
\end{equation*}%
where%
\begin{equation*}
A_{k-1}(f_{1},h):=\left(
\begin{array}{cccc}
f_{1} & \frac{1}{h+1} & \cdots & \frac{1}{(k-1)h+1} \\
\frac{1}{h+1} & \frac{1}{2h+1} & \cdots & \frac{1}{\left[ kh+1\right] } \\
\vdots & \vdots & \ddots & \vdots \\
\frac{1}{\left[ (k-1)h+1\right] } & \frac{1}{\left[ kh+1\right] } & \cdots &
\frac{1}{(2k-2)h+1}%
\end{array}%
\right)
\end{equation*}%
is now a $k\times k$ matrix. \ Continuing in this way we have%
\begin{equation*}
\begin{tabular}{l}
$\det \;A_{k}(x,h)=\left( k!(k-1)!\right) ^{2}h^{2(k+(k-1))}\left( \frac{1}{%
(k-1)h+1}\right) ^{2}\left( \frac{1}{kh+1}\right) ^{2}\left( \frac{1}{\left(
k+1\right) h+1}\right) ^{4}\cdots $ \\
\\
$\ \ \ \ \ \ \ \ \ \ \ \ \ \ \ \ \ \ \ \ \ \ \ \ \ \ \ \ \ \ \ \ \left(
\frac{1}{(2k-3)h+1}\right) ^{4}\left( \frac{1}{(2k-2)h+1}\right) ^{3}\left(
\frac{1}{(2k-1)h+1}\right) ^{2}\left( \frac{1}{2kh+1}\right) \det
\;A_{k-2}(f_{2},h)$,%
\end{tabular}%
\end{equation*}%
where
\begin{equation*}
f_{2}\equiv f_{2}(x,h,k):=f_{1}\left( \frac{(k-1)h+1}{(k-1)h}\right)
^{2}-\left( \frac{2(k-1)h+1}{(k-1)^{2}h^{2}}\right) .
\end{equation*}%
In general, we see that $\det \;A_{k}(x,h)$ can be expressed in terms of $%
\det \;A_{k-\ell -1}(f_{\ell +1},h)$, where%
\begin{equation}
f_{\ell +1}=f_{\ell }\left( \frac{(k-\ell )h+1}{\left( k-\ell \right) h}%
\right) ^{2}-\frac{2(k-\ell )h+1}{\left( \left( k-\ell \right) h\right) ^{2}}%
\text{ \ }(0\leq \ell \leq k-1).  \label{fl}
\end{equation}%
\ Thus, by direct calculation we have
\begin{equation*}
\det \;A_{k}(x,h)=\left( k^{!}\right) ^{2}h^{k(k+1)}g(h,k)f_{k},
\end{equation*}%
where $g(h,k)$ is given by (\ref{gk}) and
\begin{equation}
f_{k}=\det \;A_{0}(f_{k},h)=f_{k-1}\left( \frac{h+1}{h}\right) ^{2}-\frac{%
2h+1}{h^{2}}.  \label{fk}
\end{equation}%
On the other hand, careful inspection of the recursive definition of $f_{k}$
(cf. (\ref{fl}), (\ref{fk})) and of the formula for $f(x,h,k)$ (see (\ref%
{fx})) shows that $f_{k}=f(x,h,k)\;\;($all $x,h>0$ and $k\geq 1)$. \ The
proof is now complete.
\end{proof}

\begin{corollary}
\label{Cor1}For $k\geq 1$ and $h\geq 1$,
\begin{equation*}
\det A_{k}(x,h)<\det A_{k-1}(x,h).
\end{equation*}
\end{corollary}

\begin{proof}
We consider two cases.

\textbf{Case 1:} $\ k=1$. \ Note that $\det $ $A_{0}(x,h)=x$ and $\det $ $%
A_{1}(x,h)=\frac{x}{2h+1}-\left( \frac{1}{h+1}\right) ^{2}$. \ Thus we have%
\begin{equation*}
\frac{\det \text{ }A_{1}(x,h)}{\det \text{ }A_{0}(x,h)}<\frac{1}{2h+1}<1%
\text{.}
\end{equation*}

\textbf{Case 2:} $\ k\geq 2$. \ Consider the quotient
\begin{equation*}
\frac{\det \text{ }A_{k}(x,h)}{\det \text{ }A_{k-1}(x,h)}=\frac{%
h^{k(k+1)}\left( k^{!}\right) ^{2}g(h,k)f_{k}}{h^{k(k-1)}\left(
k-1^{!}\right) ^{2}g(h,k-1)f_{k-1}}=\frac{h^{2k}k!^{2}g(h,k)f_{k}}{%
g(h,k-1)f_{k-1}}
\end{equation*}%
and observe, using (\ref{gk}), that
\begin{eqnarray*}
\frac{g(h,k)}{g(h,k-1)} &=&\frac{\left( \frac{1}{kh+1}\right)
^{k}\prod_{j=0}^{k-1}\left( \frac{1}{(jh+1)(2kh-jh+1)}\right) ^{j+1}}{\left(
\frac{1}{\left( k-1\right) h+1}\right) ^{k-1}\prod_{j=0}^{k-2}\left( \frac{1%
}{(jh+1)(2\left( k-1\right) h-jh+1)}\right) ^{j+1}} \\
&& \\
&=&\left( \frac{1}{2kh+1}\right) \prod\nolimits_{j=0}^{k-1}\left( \frac{1}{%
\left( k+j\right) h+1}\right) ^{2}.
\end{eqnarray*}

Hence%
\begin{eqnarray*}
\frac{\det \text{ }A_{k}(x,h)}{\det \text{ }A_{k-1}(x,h)} &=&h^{2k}k!^{2}%
\left( \frac{1}{2kh+1}\right) \prod\nolimits_{j=0}^{k-1}\left( \frac{1}{%
\left( k+j\right) +1}\right) ^{2}\frac{f_{k}}{f_{k-1}} \\
&<&k!^{2}\left( \frac{1}{2kh}\right) \prod\nolimits_{j=0}^{k-1}\left( \frac{1%
}{\left( k+j\right) }\right) ^{2}\left( \frac{f_{k-1}\left( \frac{h+1}{h}%
\right) ^{2}-\frac{2h+1}{h^{2}}}{f_{k-1}}\right) \\
&<&k!^{2}\left( \frac{1}{2kh}\right) \prod\nolimits_{j=0}^{k-1}\left( \frac{1%
}{\left( k+j\right) }\right) ^{2}\left( \frac{h+1}{h}\right) ^{2} \\
&\leq &\frac{1}{4}\left( \frac{1}{2h}\right) \left( \frac{h+1}{h}\right)
^{2}<1
\end{eqnarray*}%
whenever $h\geq 1$. \ Therefore, we have $\det $ $A_{k}(x,h)<\det $ $%
A_{k-1}(x,h)\;\;$(all $h,k\geq 1$), as desired.
\end{proof}

\begin{remark}
\label{Re0} As we have mentioned before, the matrix $A_{k}(1,1)$ is the
classical Hilbert matrix. \ Specializing the above results to the case $%
x=h=1 $ in Theorem \ref{Thm-Hilbert}, we obtain%
\begin{eqnarray*}
f(1,1,k) &=&\left( k+1\right) ^{2}-\left\{ (2k-1)+(2k-3)+\cdots +(2\cdot
2+1)+(2\cdot 1+1)\right\} \\
&=&\left( k+1\right) ^{2}-k\frac{6+(k-1)2}{2}=1
\end{eqnarray*}%
and%
\begin{eqnarray*}
g(1,k) &=&\frac{0!}{\left( 2k+1\right) !}\cdot \frac{1!}{\left( 2k\right) !}%
\cdot \frac{2!}{\left( 2k-1\right) !}\cdot \frac{3!}{\left( 2k-2\right) !}%
\cdots \frac{\left( k-1\right) !}{\left( k+2\right) !}\cdot \frac{k!}{\left(
k+1\right) !} \\
&=&\frac{0!}{\left( 2k+1\right) !}\cdot \frac{1!}{\left( 2k\right) !}\cdot
\frac{2!}{\left( 2k-1\right) !}\cdot \frac{3!}{\left( 2k-2\right) !}\cdots
\frac{\left( k-1\right) !}{\left( k+2\right) !}\cdot \frac{k!}{\left(
k+1\right) !}\cdot \frac{k^{!}}{k^{!}} \\
&=&\frac{\left( k^{!}\right) ^{2}}{\left( 2k+1\right) ^{!}}.
\end{eqnarray*}%
We now use (\ref{eqHilbert}) and we recover the classical identity
\begin{equation*}
\det \;A_{k}(1,1)=\frac{\left( k^{!}\right) ^{4}}{\left( 2k+1\right) ^{!}}%
\;\;\text{(cf. \cite[Part VII, Problem 4]{PoSz}, \cite[Solution to Problem 1]%
{Choi})}.
\end{equation*}
\end{remark}

\begin{theorem}
\label{Thm(determinant)}Assume $x>0$ and $h,k\geq 1$. \ The following
statements are equivalent.\newline
(i) $\ A_{k}(x,h)\geq 0$;\newline
(ii) $\ \det \;A_{k}(x,h)\geq 0$;\newline
(iii) \ $f_{k}\equiv f(x,k,h)\geq 0$;\newline
(iv) $x\geq b(k,h)$, where
\begin{equation}
\left\{
\begin{array}{lll}
b(1,h) & := & \frac{2h+1}{\left( h+1\right) ^{2}}\text{ \ \ and} \\
b(j,h) & := & \left[ b(j-1,h)+\frac{2jh+1}{(jh)^{2}}\right] \cdot \left(
\frac{jh}{jh+1}\right) ^{2}\;\;(1\leq j\leq k).%
\end{array}%
\right.  \label{bkh}
\end{equation}
\end{theorem}

\begin{proof}
(i) $\Rightarrow $ (ii) \ This is trivial.

(ii) $\Rightarrow $ (i) \ Since
\begin{equation*}
0\leq \det \;A_{k}(x,h)<\det \;A_{k-1}(x,h)<\cdots <\det \;A_{0}(x,h)
\end{equation*}%
(by Corollary \ref{Cor1}), it follows from Choleski's Algorithm \cite{Atk}
that $A_{k}(x,h)\geq 0$.

(ii) $\Leftrightarrow $ (iii) \ This follows easily from the identity $\det
\;A_{k}(x,h)=h^{k(k+1)}\left( k^{!}\right) ^{2}g(k,h)f_{k}$ in Theorem \ref%
{Thm-Hilbert}, since $g(h,k)$ is clearly positive.

(iii) $\Leftrightarrow $ (iv) \ For $k=1$, observe that $f_{1}\geq 0\iff
\det \;A_{1}(x,h)\geq 0\iff x\equiv f_{0}\geq b(1,h)$. \ For $k\geq 2$,
recall that $f_{\ell +1}=f_{\ell }\left( \frac{(k-\ell )h+1}{\left( k-\ell
\right) h}\right) ^{2}-\frac{2(k-\ell )h+1}{\left( \left( k-\ell \right)
h\right) ^{2}}\;(0\leq \ell \leq k-1)$. \ Thus
\begin{equation*}
\begin{tabular}{l}
$f(x,k,h)\equiv f_{k}\geq 0$ \\
$\iff f_{k-1}\geq \frac{2h+1}{\left( h+1\right) ^{2}}\equiv b(1,h)$ \\
$\iff f_{k-2}\geq \left[ b(1,h)+\frac{2\left( 2h\right) +1}{\left( 2h\right)
^{2}}\right] \left( \frac{2h}{2h+1}\right) ^{2}\equiv b(2,h)$ \\
$\iff \cdots \iff f_{0}\geq \left[ b(k-1,h)+\frac{2kh+1}{(kh)^{2}}\right]
\cdot \frac{kh}{(kh+1)^{2}}\equiv b(k,h)$ \\
$\iff x\equiv f_{0}\geq b(k,h).$%
\end{tabular}%
\end{equation*}%
(Observe in passing that $b(k,h)>0$ (all $h,k\geq 1$) and that $%
\lim_{h\rightarrow \infty }b(k,h)=0\;\;$(all $k\geq 1$).)
\end{proof}

%%%%%%%%%%%%%%%%%%%%%%%%%%

\section{\label{Sect4}The class $\mathfrak{H}_{k}$ $(k\geq 1)$ is not
invariant under powers}

For a general operator $T$ on Hilbert space, it is well known that the
subnormality of $T$ implies the subnormality of $T^{m}$ $(m\geq 2)$. \ The
converse implication, however, is false; in fact, the subnormality of all
powers $T^{m}\;(m\geq 2)$ does not necessarily imply the subnormality of $T$%
, even if $T$ is a unilateral weighted shift \cite[p. 378]{Sta2}. \ Consider
for instance $W_{\alpha }\equiv \operatorname{shift}\;(a,b,1,1,\cdots )$ where $%
0<a<b<1$. \ Clearly $W_{\alpha }$ is not $2$-hyponormal (and therefore not
subnormal), but $W_{\alpha }^{m}$ is subnormal for all $m\geq 2$. \ Thus it
is indeed possible for a weighted shift $W_{\alpha }$ to have all powers $%
W_{\alpha }^{m}\;(m\geq 2)$ $k$-hyponormal without $W_{\alpha }$ being $k$%
-hyponormal. \ The example above illustrates the case $k\geq 2$. \ When $k=1$%
, it suffices to consider $W_{\alpha }\equiv \operatorname{shift}\;(1,1-x,y,y,\cdots
)$ where $0<x<1<y$. \ Then $W_{\alpha }^{m}$ $(m\geq 2)$ is hyponormal, but $%
W_{\alpha }$ is not hyponormal. \

In the multivariable case, the standard assumption on a pair $\mathbf{T}%
\equiv (T_{1},T_{2})$ is that each component $T_{i}$ is subnormal ($i=1,2$).
\ With this in mind, comparing the $k$-hyponormality of a $2$-variable
weighted shift $W_{(\alpha ,\beta )}\in \mathfrak{H}_{0}$ to the $k$%
-hyponormality of its powers $W_{(\alpha ,\beta )}^{(h,\ell )}$ is highly
nontrivial. $\ $In \cite{CLY2} we first considered this problem, for the
case of $1$-hyponormality. \ Specifically, if we let $W_{(\alpha ,\beta )}$
denote the $2$-variable weighted shift whose weight diagram is given in
Figure \ref{k-hypo}(i), we proved in \cite[Theorem 2.7]{CLY2}
that (see Figure \ref{range}): \ (i) $\ W_{(\alpha ,\beta )}$ is hyponormal 
and $W_{(\alpha ,\beta )}^{(2,1)}$ is {\textit not} hyponormal if and only if 
$a_{int}<a\leq \sqrt[4]{\frac{3}{5}}$ and $h_{21}(a)<y\leq h_{1}(a)$ ; and\
(ii) $\ W_{(\alpha ,\beta )}$ is not hyponormal and $W_{(\alpha ,\beta
)}^{(2,1)}$ is hyponormal if and only if $0<a<a_{int}$ and $h_{1}(a)<y\leq
h_{21}(a) $. \ \setlength{\unitlength}{1mm} \psset{unit=1mm}

\begin{figure}[th]
\begin{center}
\begin{picture}(110,60)
\psline{->}(0,10)(0,58)
\put(-1,10){$\_$}
\put(-1,20){$\_$}
\put(-1,30){$\_$}
\put(-1,40){$\_$}
\put(-1,50){$\_$}
\put(-6,9){$0.6$}
\put(-6,19){$0.7$}
\put(-6,29){$0.8$}
\put(-6,39){$0.9$}
\put(-6,49){$1.0$}
\put(95,6){$a$}
\put(-4,58){$\kappa$}
\psline(19.5,9.5)(19.5,10.5)
\psline(39.5,9.5)(39.5,10.5)
\psline(59.5,9.5)(59.5,10.5)
\psline(79.5,9.5)(79.5,10.5)
\put(-2,5){\footnotesize{$0$}}
\put(17.3,5){\footnotesize{$0.2$}}
\put(37.3,5){\footnotesize{$0.4$}}
\put(57.3,5){\footnotesize{$0.6$}}
\put(77.3,5){\footnotesize{$0.8$}}

 \psline{->}(83.86,3)(83.86,9)
\put(81,1){\footnotesize{$a_{int}$}}
\put(40,28){\footnotesize{$h_{1}$}}
\put(88,50){\footnotesize{$h_{1}$}}
\put(88,28){\footnotesize{$h_{21}$}}
\put(47,38){\footnotesize{$h_{21}$}}

\psline[linestyle=dashed,dash=3pt 2pt]{-}(83.86,49)(83.86,10)

\put(3,53){\footnotesize{$W_{(\alpha ,\beta )}\notin\mathfrak{H}_{1}$ and
$W_{(\alpha ,\beta )}^{(2,1)}$ hyponormal for
 $\ (a,\kappa)$ in this region}}

\psline{->}(0,10)(99,10) \psline[linestyle=dashed,dash=3pt 2pt]{-}(88.1,49)(88.1,10)
\put(50,59){\footnotesize{$W_{(\alpha ,\beta )}$ hyponormal and
$W_{(\alpha ,\beta )}^{(2,1)}\notin \mathfrak{H}_{1}$ for
 $\ (a,\kappa)$ in this region}}

\psline{->}(119,57)(87,47)

%graph of h_{1}
\pscurve[linewidth=1pt](0,23.6)(20,25.421)(40,30.50)(60,38.28)(80,47.7)(83,49)(88,49.55)

%graph of h_{21}
\pscurve[linewidth=1pt](0,25.6)(20,27.7)(40,33.1)(60,41.3)(80,49.8)(83,49.8)(85,48.06)(86,45.33)(87.79,10)

%graph of h_{12}
\pspolygon*[linecolor=lightgray](0,23.6)(20,25.421)(40,30.50)(60,38.28)(80,47.7)(83,49)(83,49.8)(80,49.8)(70,46)(60,41)(40,32.6)(20,27.3)(0,25.3)

\psline{->}(24,51)(46.0,33.5)

\end{picture}
\end{center}
\caption{Graphs of $h_{1}$ and $h_{21}$ on the interval $[0,\sqrt[4]{\frac{3%
}{5}}]$}
\label{range}
\end{figure}
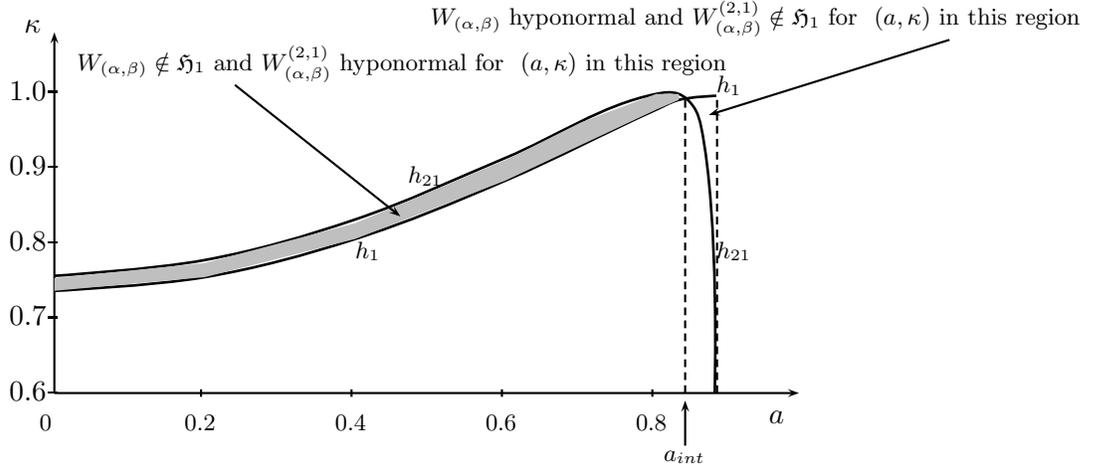
In this section we extend the above mentioned result to $k$-hyponormality (for arbitrary $k \geq 2$)
 and we also give a negative answer to Problem \ref%
{question1}. \ Our main result, Theorem \ref{Thm(invariant)-00}, gives
necessary and sufficient conditions for $W_{(\alpha ,\beta )}$ as above to
have the property $W_{(\alpha ,\beta )}\notin \mathfrak{H}_{k}$ and $%
W_{(\alpha ,\beta )}^{(h,\ell )}\in \mathfrak{H}_{k}$, for each $k\geq 2$. \

To study $k$-hyponormality of multivariable weighted shifts, we first recall
that, in one variable, the $n$-th power of a weighted shift is unitarily
equivalent to the direct sum of $n$ weighted shifts. \ Something similar
happens in two variables, as we will see in the proof of Theorem \ref%
{Thm(invariant)-00} below. \ First, we need some terminology.

Let $\mathcal{H}\equiv \ell ^{2}(\mathbb{Z}_{+})=\bigvee_{j=0}^{\infty
}\{e_{j}\}$. \ Given integers $i$ and $m$ $(m\geq 1$, $0\leq i\leq m-1)$,
define $\mathcal{H}_{i}:=\bigvee_{j=0}^{\infty }\{e_{mj+i}\}$; clearly, $%
\mathcal{H}=\bigoplus_{i=0}^{m-1}\mathcal{H}_{i}$. \ Following the notation
in \cite{CuP}, for a weight sequence $\alpha \equiv \{\alpha
_{n}\}_{n=0}^{\infty }$ we let
\begin{equation}
W_{\alpha (m:i)}:=\operatorname{shift}\;(\Pi _{n=0}^{m-1}\alpha
_{mj+i+n})_{j=0}^{\infty };  \label{ra}
\end{equation}%
that is, $W_{\alpha (m:i)}$ denotes the sequence of products of weights in
adjacent packets of size $m$, beginning with $\alpha _{i}\cdots \alpha
_{i+m-1}$. \ For example, given a weight sequence $\alpha \equiv \{\alpha
_{n}\}_{n=0}^{\infty }$, we have $W_{\alpha (2:0)}=\operatorname{shift}\;(\alpha
_{0}\alpha _{1},\alpha _{2}\alpha _{3},\cdots )$, $W_{\alpha (2:1)}=\operatorname{%
shift}\;(\alpha _{1}\alpha _{2},\alpha _{3}\alpha _{4},\cdots )$ and $%
W_{\alpha (3:2)}=\operatorname{shift}\;(\alpha _{2}\alpha _{3}\alpha _{4},\alpha
_{5}\alpha _{6}\alpha _{7},\cdots )$. \

\begin{lemma}
(\cite[Corollary 2.8]{CuP})\label{Lem1} (i) \ Let $k\geq 1$. \ Then $%
W_{\alpha }^{m}$ is $k$-hyponormal $\Leftrightarrow W_{\alpha (m:i)}$ is $k$%
-hyponormal for $0\leq i\leq m-1$.\newline
(ii) \ $W_{\alpha }^{m}$ is subnormal $\Leftrightarrow W_{\alpha (m:i)}$ is
subnormal for $0\leq i\leq m-1$.
\end{lemma}

We now introduce a key family of examples. \ Given $0<\kappa <1$, we let $%
x\equiv \{x_{n}\}_{n=0}^{\infty }$ be given by
\begin{equation}
x_{n}:=%
\begin{cases}
\kappa \sqrt{\frac{3}{4}}, & \text{if }n=0 \\
\frac{\sqrt{(n+1)(n+3)}}{(n+2)}, & \text{if }n\geq 1.%
\end{cases}
\label{alpha}
\end{equation}%
It is easy to see that $W_{x}\equiv \operatorname{shift}\;(x_{0},x_{1},x_{2},\cdots
) $ is subnormal, with Berger measure supported in $[0,1]$ and given by
\begin{equation*}
d\xi _{x}(s):=(1-\kappa ^{2})d\delta _{0}(s)+\frac{\kappa ^{2}}{2}ds+\frac{%
\kappa ^{2}}{2}d\delta _{1}(s)\;\;\text{(\cite[Proposition 4.2]{CLY1}).}
\end{equation*}%
\setlength{\unitlength}{1mm} \psset{unit=1mm}
\begin{figure}[th]
\begin{center}
\begin{picture}(165,85)

\psline{->}(20,20)(85,20)
\psline(20,40)(83,40)
\psline(20,60)(83,60)
\psline(20,80)(83,80)
\psline{->}(20,20)(20,85)
\psline(40,20)(40,83)
\psline(60,20)(60,83)
\psline(80,20)(80,83)

\put(12,16){\footnotesize{$(0,0)$}}
\put(36.5,16){\footnotesize{$(1,0)$}}
\put(56.5,16){\footnotesize{$(2,0)$}}
\put(76.5,16){\footnotesize{$(3,0)$}}

\put(27,21){\footnotesize{$x_{0}$}}
\put(47,21){\footnotesize{$x_{1}$}}
\put(67,21){\footnotesize{$x_{2}$}}
\put(81,21){\footnotesize{$\cdots$}}

\put(27,41){\footnotesize{$a$}}
\put(47,41){\footnotesize{$1$}}
\put(67,41){\footnotesize{$1$}}
\put(81,41){\footnotesize{$\cdots$}}

\put(27,61){\footnotesize{$a$}}
\put(47,61){\footnotesize{$1$}}
\put(67,61){\footnotesize{$1$}}
\put(81,61){\footnotesize{$\cdots$}}

\put(27,81){\footnotesize{$\cdots$}}
\put(47,81){\footnotesize{$\cdots$}}
\put(67,81){\footnotesize{$\cdots$}}

\psline{->}(40,10)(60,10)
\put(48,5.5){$\rm{T}_1$}
\psline{->}(10,40)(10,60)
\put(4,50){$\rm{T}_2$}

\put(11,39){\footnotesize{$(0,1)$}}
\put(11,59){\footnotesize{$(0,2)$}}
\put(11,79){\footnotesize{$(0,3)$}}

\put(20,28){\footnotesize{$y$}}
\put(20,48){\footnotesize{$1$}}
\put(20,68){\footnotesize{$1$}}
\put(21,81){\footnotesize{$\vdots$}}

\put(40,28){\footnotesize{$\frac{ay}{x_{0}}$}}
\put(40,48){\footnotesize{$1$}}
\put(40,68){\footnotesize{$1$}}
\put(41,81){\footnotesize{$\vdots$}}

\put(60,28){\footnotesize{$\frac{ay}{x_{0}x_{1}}$}}
\put(60,48){\footnotesize{$1$}}
\put(60,68){\footnotesize{$1$}}
\put(61,81){\footnotesize{$\vdots$}}

\put(30,2){$(i)$}

%\\\\\\\\\\\\\\\\\\\\\\\\\\\\\\\\\\
\put(110,2){$(ii)$}

\psline{->}(100,20)(163,20)
\psline(100,40)(162,40)
\psline(100,60)(162,60)
\psline(100,80)(162,80)
\psline{->}(100,20)(100,85)
\psline(120,20)(120,83)
\psline(140,20)(140,83)
\psline(160,20)(160,83)

\put(92,16){\footnotesize{$(0,1)$}}
\put(116,16){\footnotesize{$(1,1)$}}
\put(136,16){\footnotesize{$(2,1)$}}
\put(158,16){\footnotesize{$\cdots$}}

\put(107,21){\footnotesize{$\sqrt{\gamma_{h}}$}}
\put(127,22){\footnotesize{$\sqrt{\frac{\gamma_{2h}}{\gamma_{h}}}$}}
\put(147,22){\footnotesize{$\sqrt{\frac{\gamma_{3h}}{\gamma_{2h}}}$}}

\put(107,41){\footnotesize{$a$}}
\put(127,41){\footnotesize{$1$}}
\put(147,41){\footnotesize{$1$}}

\put(107,61){\footnotesize{$a$}}
\put(127,61){\footnotesize{$1$}}
\put(147,61){\footnotesize{$1$}}

\put(106,81){\footnotesize{$\cdots$}}
\put(127,81){\footnotesize{$\cdots$}}
\put(147,81){\footnotesize{$\cdots$}}

\psline{->}(120,10)(140,10)
\put(126,6){$T_{1}^{h}|_{_{\mathcal{H}_{0,0}}}$}

\psline{->}(96.5,40)(96.5,60)
\put(85,50){$T_{2}^{\ell}|_{_{\mathcal{H}_{0,0}}}$}

\put(100,29){\footnotesize{$y$}}
\put(100,49){\footnotesize{$1$}}
\put(100,69){\footnotesize{$1$}}
\put(101,83){\footnotesize{$\vdots$}}

\put(120,29){\footnotesize{$\frac{ay}{\sqrt{\gamma_{h}}}$}}
\put(120,48){\footnotesize{$1$}}
\put(120,68){\footnotesize{$1$}}
\put(121,83){\footnotesize{$\vdots$}}

\put(140,29){\footnotesize{$\frac{ay}{\sqrt{\gamma_{2h}}}$}}
\put(140,48){\footnotesize{$1$}}
\put(140,68){\footnotesize{$1$}}
\put(141,83){\footnotesize{$\vdots$}}

\end{picture}
\end{center}
\caption{Weight diagram of the $2$-variable weighted shifts in Theorem
\protect\ref{examplethm}, \protect\ref{Thm-main}, \protect\ref{Thm8} and
weight diagram of the $2 $-variable weighted shift $W_{(\protect\alpha ,%
\protect\beta )}^{(h,\ell )}$ , respectively.}
\label{k-hypo}
\end{figure}
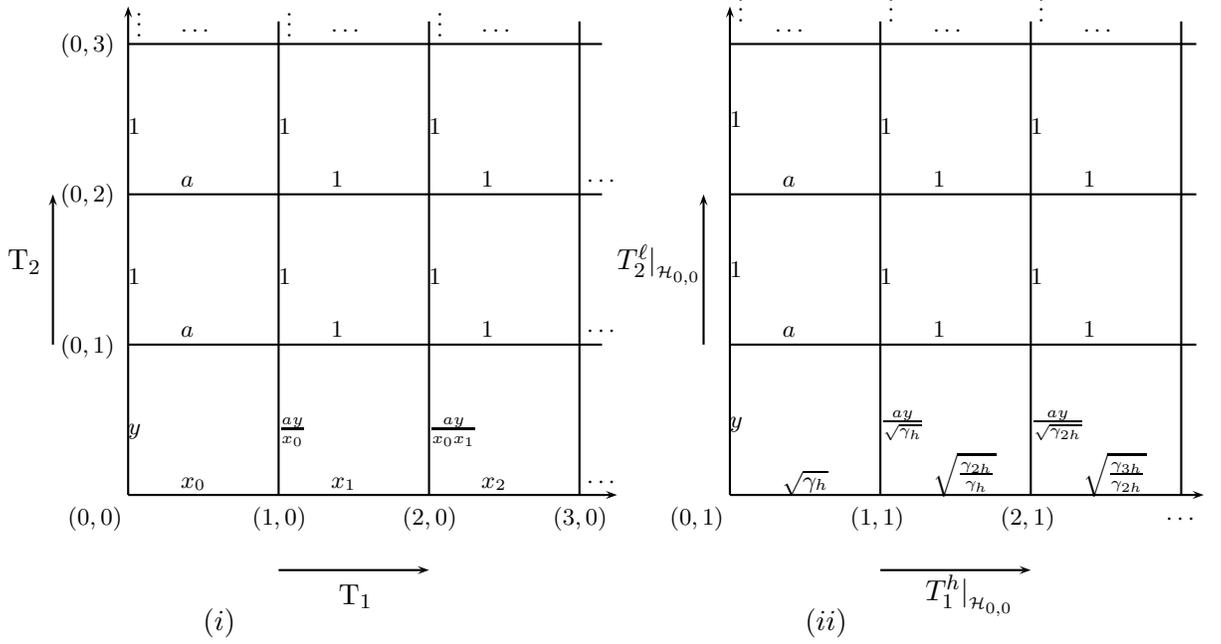
Consider now the $2$-variable weighted shift given in Figure \ref{k-hypo}%
(i), where $W_{x}\equiv \operatorname{shift}\;(x_{0},x_{1},x_{2},\cdots )$ and $%
y:=\kappa $.

\begin{theorem}
(\cite{CLY1})\label{examplethm} \ For $0<a\leq \frac{1}{\sqrt{2}}$, $%
0<\kappa <1$, $x_{n}$ as in (\ref{alpha}) and $y:=\kappa $, let $W_{(\alpha
,\beta )}\equiv (T_{1},T_{2})$ be the $2$-variable weighted shift given by
Figure \ref{k-hypo}(i). \ For $k\geq 2$, let
\begin{equation}
\begin{tabular}{l}
$F(a,k):=\sqrt{\frac{\frac{(k+1)^{2}}{2k(k+2)}-a^{2}}{a^{4}-\frac{5}{2}a^{2}+%
\frac{(k+1)^{2}}{2k(k+2)}+\frac{2k^{2}+4k+3}{4(k+1)^{2}}}}.$%
\end{tabular}
\label{fak}
\end{equation}%
Then\newline
(i) $\ T_{1}$ and $T_{2}$ are subnormal;\newline
(ii) \ $W_{(\alpha ,\beta )}\in \mathfrak{H}_{1}\Leftrightarrow 0<\kappa
\leq h_{1}(a):=\sqrt{\frac{32-48a^{4}}{59-72a^{2}}}$;\newline
(iii) $\ W_{(\alpha ,\beta )}\in \mathfrak{H}_{k}\Leftrightarrow 0<\kappa
\leq h_{k}(a):=F(a,k)\;(k\geq 2)$;\newline
(iv) $\ W_{(\alpha ,\beta )}\in \mathfrak{H}_{\infty }\Leftrightarrow
0<\kappa \leq h_{\infty}(a):=\sqrt{\frac{1}{2-a^{2}}}$.\newline
In particular, $W_{(\alpha ,\beta )}$ is hyponormal and not subnormal if and
only if $\sqrt{\frac{1}{2-a^{2}}}<\kappa \leq \sqrt{\frac{32-48a^{4}}{%
59-72a^{2}}}$.
\end{theorem}

We now recall that, by (\ref{bkh}),
\begin{equation*}
b(k,h)=\left[ b(k-1,h)+\frac{2kh+1}{(kh)^{2}}\right] \cdot \left( \frac{kh}{%
kh+1}\right) ^{2}\text{ and }b(1,h)=\frac{2h+1}{(h+1)^{2}}.
\end{equation*}%
Using mathematical induction we can see that
\begin{equation}
b(k,1)=\frac{k(k+2)}{\left( k+1\right) ^{2}}.  \label{bk1}
\end{equation}

\begin{remark}
\label{Re-00}If $x=\frac{2(1-\kappa ^{2})}{\kappa ^{2}}$ in Theorem \ref%
{Thm(determinant)}, then for $h\geq 1$ and $k\geq 1$,%
\begin{equation*}
\det \;A_{k}(x,h)\geq 0\Leftrightarrow x\geq b(k,h)\Leftrightarrow \kappa
\leq \sqrt{\frac{2}{2+b(k,h)}}.
\end{equation*}
\end{remark}

\begin{lemma}
\label{Le-00}For $h\geq 1$ and $k\geq 1$, we have $b(k,h)\leq b(k,1)$.
\end{lemma}

\begin{proof}
We fix $h\geq 1$ and use induction on $k$. \ First, observe that, on the
interval $[1,+\infty ),$ $b(1,h)\equiv \frac{2h+1}{(h+1)^{2}}$ is a
decreasing function of $h$, so we clearly have $b(1,h)\leq b(1,1)$. \ For
the induction step, assume now that $k\geq 2$ and that $b(k-1,h)\leq
b(k-1,1) $. \ Then
\begin{eqnarray*}
b(k,h) &=&\left[ b(k-1,h)+\frac{2kh+1}{(kh)^{2}}\right] \cdot \left( \frac{kh%
}{kh+1}\right) ^{2} \\
&\leq &\left[ b(k-1,1)+\frac{2kh+1}{(kh)^{2}}\right] \cdot \left( \frac{kh}{%
kh+1}\right) ^{2} \\
&=&\left[ \frac{(k-1)(k+1)}{k^{2}}+\frac{2kh+1}{(kh)^{2}}\right] \cdot
\left( \frac{kh}{kh+1}\right) ^{2}\;\;\text{(by (\ref{bk1}))} \\
&=&\frac{(kh+1)^{2}-h^{2}}{(kh+1)^{2}}\leq \frac{k(k+2)}{(k+1)^{2}},
\end{eqnarray*}%
since the next-to-the-last expression is a decreasing function of $h$ on the
interval $[1,+\infty )$. \ We therefore have $b(k,h)\leq \frac{k(k+2)}{%
(k+1)^{2}}\equiv b(k,1)$, as desired.
\end{proof}

\begin{corollary}
\label{Cor-01}For $h\geq 1$ and $k\geq 1$,
\begin{equation}
F\left( \sqrt{\tfrac{1}{2}},k\right) =\sqrt{\frac{2}{2+b(k,1)}}\leq \sqrt{%
\frac{2}{2+b(k,h)}}.  \label{cor45}
\end{equation}
\end{corollary}

\begin{proof}
From Lemma \ref{Le-00} we know that $b(k,h)\leq b(k,1)$. \ Thus it suffices
to establish in (\ref{cor45}). \ A direct calculation using (\ref{fak})
shows that $F\left( \sqrt{\frac{1}{2}},k\right) ^{2}=\frac{2(k+1)^{2}}{%
3k^{2}+6k+2}$, and from (\ref{bk1}) we know that $b(k,1)=\frac{k(k+2)}{%
(k+1)^{2}}$. \ It follows that
\begin{eqnarray*}
\frac{2}{2+b(k,1)}-F\left( \sqrt{\tfrac{1}{2}},k\right) ^{2} &=&\frac{2}{2+%
\frac{k(k+2)}{(k+1)^{2}}}-\frac{2(k+1)^{2}}{3k^{2}+6k+2} \\
&=&\frac{2(k+1)^{2}}{2(k+1)^{2}+k(k+2)}-\frac{2(k+1)^{2}}{3k^{2}+6k+2}=0,
\end{eqnarray*}%
as desired.
\end{proof}

\begin{remark}
\label{Re-01}From Lemma \ref{Le-00} and Corollary \ref{Cor-01} we see at
once that for $h\geq 1$ and $k\geq 1$, the following statements are
equivalent:\newline
(i) $\ b(k,h)<b(k,1)$;\newline
(ii) $\ F\left( \sqrt{\frac{1}{2}},k\right) <\sqrt{\frac{2}{2+b(k,h)}}$.
\end{remark}

\begin{lemma}
\label{diff}Let $G(h):=\frac{2h^{3}+7h^{2}+8h+3}{2h^{3}+7h^{2}+10h+4}$. \
Then $G$ is an increasing function of $h$ on $[1,\infty )$, $G(1)=\frac{20}{%
23}$ and $\lim_{h\rightarrow \infty }G(h)=1$.
\end{lemma}

\begin{proof}
$\lim_{h\rightarrow \infty }G(h)=1$ is clear. \ To establish that $G$ is
increasing, observe that%
\begin{equation*}
G^{\prime }(h)=\frac{2\left( 4h^{3}+10h^{2}+7h+1\right) }{\left(
2h^{3}+7h^{2}+10h+4\right) ^{2}}>0
\end{equation*}%
on $[1,\infty )$.
\end{proof}

We are now ready to prove our main result of this section.

\begin{theorem}
\label{Thm(invariant)-00}Let $W_{(\alpha ,\beta )}\equiv (T_{1},T_{2})$ be
the $2$-variable weighted shift whose weight diagram is given in Figure \ref%
{k-hypo}(i) (where $a=\sqrt{\frac{1}{2}}$ and $W_{x}$ is as in (\ref{alpha}%
)). \ Then given $k,\ell \geq 1$ and $h\geq 2$,%
\begin{equation*}
W_{(\alpha ,\beta )}^{(h,\ell )}\in \mathfrak{H}_{k}\text{ but }W_{(\alpha
,\beta )}\notin \mathfrak{H}_{k}\Leftrightarrow \left\{
\begin{tabular}{ll}
$\sqrt{\frac{20}{23}}<\kappa \leq \sqrt{\frac{63}{68}}$, & if $k=1$ \\
$\sqrt{\frac{2(k+1)^{2}}{3k^{2}+6k+2}}<\kappa \leq \sqrt{\frac{2}{2+b(k,h)}}$%
, & if $k\geq 2.$%
\end{tabular}%
\right.
\end{equation*}
\end{theorem}

\begin{proof}
From Lemma \ref{khypo}, we recall that a $2$-variable weighted shift $%
W_{(\alpha ,\beta )}$ is $k$-hyponormal if and only if
\begin{equation}
M_{\mathbf{k}}(k)=(\gamma _{\mathbf{k}+(m,n)+(p,q)})_{_{0\leq p+q\leq
k}^{0\leq n+m\leq k}}\geq 0,  \label{k-hy}
\end{equation}%
for all $\mathbf{k}\in \mathbb{Z}_{+}^{2}$.

We first let $\mathcal{H}_{(m,n)}:=\bigvee_{i,j=0}^{\infty }\{e_{(hi+m,\ell
j+n)}:h,\ell \geq 1\}$, for $0\leq m\leq h-1$ and $0\leq n\leq \ell -1$. \
Then we have $\ell ^{2}(\mathbb{Z}_{+}^{2})\equiv
\bigoplus_{m=0}^{h-1}\bigoplus_{n=0}^{\ell -1}\mathcal{H}_{(m,n)}$. $\ $%
Observe that $\mathcal{H}_{(m,n)}$ reduces $T_{1}^{h}$ and $T_{2}^{\ell }$.
\ Thus if a $2$-variable weighted shift $W_{(\alpha ,\beta )}$ is given as
in Figure \ref{k-hypo}(i), then for $h,\ell \geq 1$, we can write%
\begin{equation*}
W_{(\alpha ,\beta )}^{(h,\ell )}\equiv (T_{1}^{h},T_{2}^{\ell })\cong
(W_{\alpha (h:0)}\oplus (I\otimes S_{\sqrt{\frac{1}{2}}}),T_{2}|_{\mathcal{H}%
_{0}})\bigoplus \bigoplus_{i=1}^{h-1}(W_{\alpha (h:i)}\oplus (I\otimes
U_{+}),T_{2}|_{\mathcal{H}_{i}}),
\end{equation*}%
where%
\begin{equation*}
W_{\alpha (h:i)}=\operatorname{shift}\;(\sqrt{\frac{\gamma _{(i+1)h}}{\gamma _{ih}}},%
\sqrt{\frac{\gamma _{(i+2)h}}{\gamma _{(i+1)h}}},\cdots )\text{ and }%
\mathcal{H}_{i}:=\bigoplus_{n=0}^{\ell -1}\mathcal{H}_{(i,n)}\;\;\text{(}%
0\leq i\leq k-1\text{)}.
\end{equation*}%
Clearly, $\left\Vert W_{\alpha (h:i)}\right\Vert =1$, and the Berger measure
of $W_{\alpha (h:i)}$ has an atom at $1$, so by Lemma \ref{thmbackext} we
see that $(W_{\alpha (h:i)}\oplus (I\otimes U_{+}),T_{2}|_{\mathcal{H}_{i}})$
$(1\leq i\leq h-1)$ is subnormal. \ Thus, for $k\geq 1$, the $k$%
-hyponormality of $W_{(\alpha ,\beta )}^{(h,\ell )}$ is equivalent to the $k$%
-hyponormality of $(W_{\alpha (h:0)}\oplus (I\otimes S_{\sqrt{\frac{1}{2}}%
}),T_{2}|_{\mathcal{H}_{0}}).$ \ Observe that%
\begin{equation*}
\begin{tabular}{l}
$(W_{\alpha (h:0)}\oplus (I\otimes S_{\sqrt{\frac{1}{2}}}),T_{2}|_{\mathcal{H%
}_{0}})\cong \bigoplus_{n=0}^{\ell -1}(W_{\alpha (h:0)}\oplus (I\otimes S_{%
\sqrt{\frac{1}{2}}}),T_{2}^{\ell }|_{\mathcal{H}_{(0,n)}})$%
\end{tabular}%
\end{equation*}%
and%
\begin{equation*}
\begin{tabular}{l}
$\bigoplus_{n=0}^{\ell -1}(W_{\alpha (h:0)}\oplus (I\otimes S_{\sqrt{\frac{1%
}{2}}}),T_{2}^{\ell }|_{\mathcal{H}_{(0,n)}})\cong (W_{\alpha (h:0)}\oplus
(I\otimes S_{\sqrt{\frac{1}{2}}}),T_{2}^{\ell }|_{\mathcal{H}%
_{(0,0)}})\bigoplus \bigoplus_{n=0}^{\ell -1}(I\otimes S_{\sqrt{\frac{1}{2}}%
},U_{+}\otimes I)$.%
\end{tabular}%
\end{equation*}%
Observe that the second summand is clearly subnormal; thus, for $h,\ell \geq
1 $, the $k$-hyponormality of $(T_{1}^{h},T_{2}^{\ell })$ is equivalent to
the $k$-hyponormality of the first summand, $(W_{\alpha (h:0)}\oplus
(I\otimes S_{\sqrt{\frac{1}{2}}}),T_{2}^{\ell }|_{\mathcal{H}_{(0,0)}})$. \
Observe also that%
\begin{equation*}
(W_{\alpha (h:0)}\oplus (I\otimes S_{\sqrt{\frac{1}{2}}}),T_{2}^{\ell }|_{%
\mathcal{H}_{(0,0)}})\cong (W_{\alpha (h:0)}\oplus (I\otimes S_{\sqrt{\frac{1%
}{2}}}),T_{2}|_{\mathcal{H}_{(0,0)}}).
\end{equation*}%
Thus
\begin{equation*}
W_{(\alpha ,\beta )}^{(h,\ell )}\in \mathfrak{H}_{k}\Leftrightarrow
(W_{\alpha (h:0)}\oplus (I\otimes S_{\sqrt{\frac{1}{2}}}),T_{2}|_{\mathcal{H}%
_{(0,0)}})\in \mathfrak{H}_{k}.
\end{equation*}%
We consider two cases.

\textbf{Case 1:} $k=1$. \ To check hyponormality, by Lemma \ref{joint hypo}
and Lemma \ref{thmbackext} it suffices to apply the Six-point Test at $%
\mathbf{k=(}0,0)$. \ A direct calculation shows that%
\begin{equation*}
\begin{tabular}{l}
$H_{(W_{\alpha (h:0)}\oplus (I\otimes S_{\sqrt{\frac{1}{2}}}),T_{2}|_{%
\mathcal{H}_{(0,0)}})}(0,0)\geq 0\Leftrightarrow \kappa \leq G(h)\equiv
\sqrt{\frac{2h^{3}+7h^{2}+8h+3}{2h^{3}+7h^{2}+10h+4}}$%
\end{tabular}%
\end{equation*}%
(cf. Lemma \ref{diff}). \ Therefore, for all $h,\ell \geq 1$, we have%
\begin{equation*}
W_{(\alpha ,\beta )}^{(h,\ell )}\in \mathfrak{H}_{1}\Leftrightarrow \kappa
\leq G(h).
\end{equation*}%
Since $G(h)$ is an increasing function, we see that if $\sqrt{\frac{20}{23}}%
=G(1)<\kappa \leq G(2)=\sqrt{\frac{63}{68}}$, we simultaneously get $%
W_{(\alpha ,\beta )}^{(h,\ell )}\in \mathfrak{H}_{1}$ and $W_{(\alpha ,\beta
)}\notin \mathfrak{H}_{1}$ (all $h\geq 2,\ell \geq 1$).

\textbf{Case 2:} $k\geq 2$. \ Note that%
\begin{equation*}
W_{(\alpha ,\beta )}^{(h,\ell )}\in \mathfrak{H}_{k}\Leftrightarrow
(W_{\alpha (h:0)}\oplus (I\otimes S_{\sqrt{\frac{1}{2}}}),T_{2}|_{\mathcal{H}%
_{(0,0)}})\in \mathfrak{H}_{k}.
\end{equation*}%
To check the $k$-hyponormality of $(W_{\alpha (h:0)}\oplus (I\otimes S_{%
\sqrt{\frac{1}{2}}}),T_{2}|_{\mathcal{H}_{(0,0)}})$, we observe that it
suffices to apply Lemma \ref{khypo}(ii) at $\mathbf{k=(}0,0)$. \ Now, the
moments associated with $(W_{\alpha (h:0)}\oplus (I\otimes S_{\sqrt{\frac{1}{%
2}}}),T_{2}|_{\mathcal{H}_{(0,0)}})$ are
\begin{equation}
\gamma _{\mathbf{k}}(W_{(\alpha ,\beta )}^{(h,\ell )})=\left\{
\begin{tabular}{ll}
$1,$ & $\text{if }k_{1}=0\text{ and }k_{2}=0$ \\
$\gamma _{k_{1}h}(W_{(\alpha ,\beta )}),$ & $\text{if }k_{1}\geq 1\text{ and
}k_{2}=0$ \\
$\kappa ^{2},$ & $\text{if }k_{1}=0\text{ and }k_{2}\geq 1$ \\
$\frac{\kappa ^{2}}{2},$ & $\text{if }k_{1}\geq 1\text{ and }k_{2}\geq 1$.%
\end{tabular}%
\right.  \label{moment0}
\end{equation}%
By direct computation (i.e., interchanging rows and columns, discarding some
redundant rows and columns, and multiplying by $\frac{2}{\kappa ^{2}}$ in
the moment matrix of $(W_{\alpha (h:0)}\oplus (I\otimes S_{\sqrt{\frac{1}{2}}%
}),T_{2}|_{\mathcal{H}_{(0,0)}})$), we see that for $0<\kappa <1$ and $%
h,\ell \geq 1$,%
\begin{equation*}
\begin{tabular}{l}
$W_{(\alpha ,\beta )}^{(h,\ell )}\in \mathfrak{H}_{k}$ \\
$\Leftrightarrow W_{(\alpha ,\beta )}^{(h,\ell )}|_{\mathcal{H}_{(0,0)}}\in
\mathfrak{H}_{k}$ \\
$\Leftrightarrow (W_{\alpha (h:0)}\oplus (I\otimes S_{\sqrt{\frac{1}{2}}%
}),T_{2}|_{\mathcal{H}_{(0,0)}})\in \mathfrak{H}_{k}$ \\
$\Leftrightarrow J_{k}(\kappa ,h)\geq 0$ \\
$\Leftrightarrow L_{k}(\kappa ,h)\geq 0$,%
\end{tabular}%
\end{equation*}%
where%
\begin{equation*}
\begin{tabular}{l}
$J_{k}(\kappa ,h):=$ \\
\\
$\left(
\begin{array}{cccccccc}
1 & \frac{\kappa ^{2}(h+2)}{2(h+1)} & \kappa ^{2} & \frac{\kappa ^{2}(2h+2)}{%
2(2h+1)} & \frac{\kappa ^{2}}{2} & \frac{\kappa ^{2}(3h+2)}{2(3h+1)} & \cdots
& \frac{\kappa ^{2}(kh+2)}{2(kh+1)} \\
\frac{\kappa ^{2}(h+2)}{2(h+1)} & \frac{\kappa ^{2}(2h+2)}{2(2h+1)} & \frac{%
\kappa ^{2}}{2} & \frac{\kappa ^{2}(3h+2)}{2(3h+1)} & \frac{\kappa ^{2}}{2}
& \frac{\kappa ^{2}(4h+2)}{2(4h+1)} & \cdots & \frac{\kappa ^{2}((k+1)h+2)}{%
2((k+1)h+1)} \\
\kappa ^{2} & \frac{\kappa ^{2}}{2} & \kappa ^{2} & \frac{\kappa ^{2}}{2} &
\frac{\kappa ^{2}}{2} & \frac{\kappa ^{2}}{2} & \cdots & \frac{\kappa ^{2}}{2%
} \\
\frac{\kappa ^{2}(2h+2)}{2(2h+1)} & \frac{\kappa ^{2}(3h+2)}{2(3h+1)} &
\frac{\kappa ^{2}}{2} & \frac{\kappa ^{2}(4h+2)}{2(4h+1)} & \frac{\kappa ^{2}%
}{2} & \frac{\kappa ^{2}(5h+2)}{2(5h+1)} & \cdots & \frac{\kappa
^{2}((k+2)h+2)}{2((k+2)h+1)} \\
\frac{\kappa ^{2}}{2} & \frac{\kappa ^{2}}{2} & \frac{\kappa ^{2}}{2} &
\frac{\kappa ^{2}}{2} & \frac{\kappa ^{2}}{2} & \frac{\kappa ^{2}}{2} &
\cdots & \frac{\kappa ^{2}}{2} \\
\frac{\kappa ^{2}(3h+2)}{2(3h+1)} & \frac{\kappa ^{2}(4h+2)}{2(4h+1)} &
\frac{\kappa ^{2}}{2} & \frac{\kappa ^{2}(5h+2)}{2(5h+1)} & \frac{\kappa ^{2}%
}{2} & \frac{\kappa ^{2}(6h+2)}{2(6h+1)} & \cdots & \frac{\kappa
^{2}((k+3)h+2)}{2((k+3)h+1)} \\
\vdots & \vdots & \vdots & \vdots & \vdots & \vdots & \ddots & \vdots \\
\frac{y^{2}(kh+2)}{2(kh+1)} & \frac{y^{2}((k+1)h+2)}{2((k+1)h+1)} & \frac{%
\kappa ^{2}}{2} & \frac{y^{2}((k+2)h+2)}{2((k+2)h+1)} & \frac{\kappa ^{2}}{2}
& \frac{\kappa ^{2}((k+3)h+2)}{2((k+3)h+1)} & \cdots & \frac{\kappa
^{2}(2kh+2)}{2(2kh+1)}%
\end{array}%
\right) _{(k+3)\times (k+3)}$%
\end{tabular}%
\end{equation*}%
and%
\begin{equation*}
\begin{tabular}{l}
$L_{k}(\kappa ,h):=$ \\
\\
$\left(
\begin{array}{cc}
\left(
\begin{array}{cc}
1 & 1 \\
1 & 2%
\end{array}%
\right) & \left(
\begin{array}{ccccc}
1 & 1 & \cdots & 1 & 1 \\
2 & 1 & \cdots & 1 & 1%
\end{array}%
\right) \\
\left(
\begin{array}{cc}
1 & 2 \\
1 & 1 \\
\vdots & \vdots \\
1 & 1 \\
1 & 1%
\end{array}%
\right) & \left(
\begin{array}{ccccc}
\frac{2}{\kappa ^{2}} & \frac{1}{h+1}+1 & \cdots & \frac{1}{(k-1)h+1}+1 &
\frac{1}{kh+1}+1 \\
\frac{1}{h+1}+1 & \frac{1}{2h+1}+1 & \cdots & \frac{1}{kh+1}+1 & \frac{1}{%
(k+1)h+1}+1 \\
\vdots & \vdots & \ddots & \vdots & \vdots \\
\frac{1}{(k-1)h+1}+1 & \frac{1}{kh+1}+1 & \cdots & \frac{1}{(2k-2)h+1}+1 &
\frac{1}{(2k-1)h+1}+1 \\
\frac{1}{kh+1}+1 & \frac{1}{(k+1)h+1}+1 & \cdots & \frac{1}{(2k-1)h+1}+1 &
\frac{1}{2kh+1}+1%
\end{array}%
\right)%
\end{array}%
\right) _{_{(k+3)\times (k+3)}}$%
\end{tabular}%
.
\end{equation*}%
Note that $\det \;\left(
\begin{array}{cc}
1 & 1 \\
1 & 2%
\end{array}%
\right) >0$, and let
\begin{equation*}
M_{k}(\kappa ,h):=\left(
\begin{array}{ccccc}
\frac{2}{\kappa ^{2}} & \frac{1}{h+1}+1 & \cdots & \frac{1}{(k-1)h+1}+1 &
\frac{1}{kh+1}+1 \\
\frac{1}{h+1}+1 & \frac{1}{2h+1}+1 & \cdots & \frac{1}{kh+1}+1 & \frac{1}{%
(k+1)h+1}+1 \\
\vdots & \vdots & \ddots & \vdots & \vdots \\
\frac{1}{(k-1)h+1}+1 & \frac{1}{kh+1}+1 & \cdots & \frac{1}{(2k-2)h+1}+1 &
\frac{1}{(2k-1)h+1}+1 \\
\frac{1}{kh+1}+1 & \frac{1}{(k+1)h+1}+1 & \cdots & \frac{1}{(2k-1)h+1}+1 &
\frac{1}{2kh+1}+1%
\end{array}%
\right) _{_{(k+1)\times (k+1)}}.
\end{equation*}%
Then we have%
\begin{equation*}
M_{k}(\kappa ,h)-\left(
\begin{array}{ccccc}
1 & 1 & \cdots & 1 & 1 \\
2 & 1 & \cdots & 1 & 1%
\end{array}%
\right) _{_{(k+1)\times (k+1)}}^{\ast }\left(
\begin{array}{cc}
1 & 1 \\
1 & 2%
\end{array}%
\right) ^{-1}\left(
\begin{array}{ccccc}
1 & 1 & \cdots & 1 & 1 \\
2 & 1 & \cdots & 1 & 1%
\end{array}%
\right) _{_{(k+1)\times (k+1)}}=:A_{k}(x,h)\text{.}
\end{equation*}%
where $x:=\frac{2(1-\kappa ^{2})}{\kappa ^{2}}$ and $A_{k}(x,h)$ is as in
Theorem \ref{Thm-Hilbert}. \ Thus, after we apply Smul'jan Lemma (Lemma \ref%
{smu}) to $L_{k}(\kappa ,h)$, we show that for $0<\kappa <1$ and $h,\ell
\geq 1$, $L_{k}(\kappa ,h)\geq 0\Leftrightarrow A_{k}(x,h)\geq 0$. \
Therefore,
\begin{equation}
\begin{tabular}{l}
$W_{(\alpha ,\beta )}^{(h,\ell )}\in \mathfrak{H}_{k}\Leftrightarrow
W_{(\alpha ,\beta )}^{(h,\ell )}|_{\mathcal{H}_{(0,0)}}\in \mathfrak{H}%
_{k}\Leftrightarrow L_{k}(\kappa ,h)\geq 0\Leftrightarrow A_{k}(x,h)\geq 0.$%
\end{tabular}
\label{iff}
\end{equation}%
From Remark \ref{Re-01}(ii), for $k,h\geq 2$, we recall that%
\begin{equation}
\begin{tabular}{l}
$b(k,h)<b(k,1)\Leftrightarrow F\left( \sqrt{\frac{1}{2}},k\right) =\sqrt{%
\frac{2(k+1)^{2}}{3k^{2}+6k+2}}<\sqrt{\frac{2}{2+b(k,2)}}$.%
\end{tabular}
\label{i-0}
\end{equation}%
By Theorem \ref{examplethm}(iii),
\begin{equation}
W_{(\alpha ,\beta )}\in \mathfrak{H}_{k}\Leftrightarrow 0<\kappa \leq
F\left( \sqrt{\tfrac{1}{2}},k\right) \;(k\geq 2)\text{.}  \label{i-1}
\end{equation}%
Now, Remark \ref{Re-00}(i) and Theorem \ref{Thm(determinant)} imply that for
$h\geq 1$ and $k\geq 2$,%
\begin{equation}
\det \;A_{k}(x,h)\geq 0\Leftrightarrow A_{k}(x,h)\geq 0\Leftrightarrow x\geq
b(k,h)\Leftrightarrow \kappa \leq \sqrt{\frac{2}{2+b(k,h)}}.  \label{i-2}
\end{equation}%
Thus, by (\ref{iff}) and (\ref{i-2}), we have that for $h,\ell \geq 1$,%
\begin{equation}
W_{(\alpha ,\beta )}^{(h,\ell )}\in \mathfrak{H}_{k}\Leftrightarrow 0<\kappa
\leq \sqrt{\frac{2}{2+b(k,h)}}\;\;(k\geq 2).  \label{i-3}
\end{equation}%
Therefore, by (\ref{iff}), (\ref{i-0}), (\ref{i-1}), (\ref{i-2}) and (\ref%
{i-3}), for $\ell \geq 1$ and $h,k\geq 2$, we have
\begin{equation}
\begin{tabular}{l}
$W_{(\alpha ,\beta )}^{(h,\ell )}\in \mathfrak{H}_{k}$ but $W_{(\alpha
,\beta )}\notin \mathfrak{H}_{k}$ if and only if $F\left( \sqrt{\frac{1}{2}}%
,k\right) <\kappa \leq \sqrt{\frac{2}{2+b(k,h)}}$,%
\end{tabular}
\label{i-4}
\end{equation}%
as desired.
\end{proof}

\begin{remark}
\label{Re-010}(i) \ We know that for $k,h\geq 1$, $b(k,h)\leq b(k,1)=\frac{%
k(k+2)}{(k+1)^{2}}$, so that $\lim \;\sup {}_{k}\;b(k,h)\leq 1$. \ As an
application of (\ref{i-3}) we can establish that $\lim_{k}\;b(k,h)$ exists.
\ Recall that $W_{(\alpha ,\beta )}\in \mathfrak{H}_{k+1}\Rightarrow
W_{(\alpha ,\beta )}\in \mathfrak{H}_{k}$, so that from (\ref{i-3}) we see
that for each fixed $h\geq 1$, $b(k,h)$ must be a nondecreasing function of $%
k$, and therefore $b(h):=\lim_{k}\;b(k,h)=\lim \;\sup {}_{k}\;b(k,h)\leq 1$.%
\newline
(ii) \ We believe it is nontrivial to show that for $h\geq 1$, $%
\lim_{k\rightarrow \infty }$ $b(k,h)=1$. \ We now provide an
operator-theoretic proof of this fact. \ By (\ref{i-1}) and (\ref{i-3}), for
$k\geq 2$, we have $W_{(\alpha ,\beta )}^{(h,\ell )}\in \mathfrak{H}%
_{k}\Leftrightarrow 0<\kappa \leq \sqrt{\frac{2}{2+b(k,h)}}$. \ Since $%
b(k,h) $ is a nondecreasing function of $k$, and $\lim {}_{k\rightarrow
\infty }\;b(k,h)=b(h)\leq 1$, we easily see that
\begin{equation}
W_{(\alpha ,\beta )}^{(h,\ell )}\in \mathfrak{H}_{\infty }\iff 0<\kappa \leq
\sqrt{\frac{2}{2+b(h)}}.  \label{eq11}
\end{equation}%
We now let $\mathcal{M}_{1}(0,0)$ denote the subspace of $\mathcal{H}%
_{(0,0)} $ spanned by canonical orthonormal basis vectors with indices $%
\mathbf{k}=(k_{1},k_{2})$ with $k_{1}\geq 0$ and $k_{2}\geq 1$. \ We have $%
W_{(\alpha ,\beta )}^{(h,\ell )}\mathbf{|}_{\mathcal{M}_{1}(0,0)}\cong
(I\otimes S_{a},U_{+}\otimes I)\in \mathfrak{H}_{\infty }$ with Berger
measure $\mu _{\mathcal{M}_{1}(0,0)}:=[(1-a^{2})\delta _{0}+a^{2}\delta
_{1}]\times \delta _{1}$. \ Thus, by Lemma \ref{backext} and a direct
calculation, we see that%
\begin{equation*}
W_{(\alpha ,\beta )}^{(h,\ell )}\in \mathfrak{H}_{\infty }\Leftrightarrow
W_{(\alpha ,\beta )}^{(h,\ell )}\mathbf{|}_{\mathcal{H}(0,0)}\in \mathfrak{H}%
_{\infty }\Leftrightarrow \kappa ^{2}\left\Vert \frac{1}{t}\right\Vert
_{L^{1}(\mu _{\mathcal{M}_{1}(0,0)})}(\mu _{\mathcal{M}_{1}(0,0)})_{ext}^{X}%
\leq \xi _{x}\Leftrightarrow 0<\kappa \leq \sqrt{\tfrac{2}{3}}\text{,}
\end{equation*}%
that is, for $h,\ell \geq 1$,
\begin{equation}
W_{(\alpha ,\beta )}^{(h,\ell )}\in \mathfrak{H}_{\infty }\Leftrightarrow
0<\kappa \leq \sqrt{\tfrac{2}{3}}.  \label{eq12}
\end{equation}%
From (\ref{eq11}) and (\ref{eq12}) we see at once that $b(h)=1$, as desired.
\end{remark}

\begin{example}
\label{Re-02}As specific instances of Theorem \ref{Thm(invariant)-00}, we
have\newline
(i) $\ $%
\begin{eqnarray*}
(W_{(\alpha ,\beta )}^{(9,1)} &\in &\mathfrak{H}_{2}\text{ and }W_{(\alpha
,\beta )}\notin \mathfrak{H}_{1}) \\
&\Leftrightarrow &0.932505\simeq \sqrt{\tfrac{20}{23}}<\kappa \leq \sqrt{%
\frac{2}{2+b(2,9)}}=\sqrt{\tfrac{9025}{10257}}\simeq 0.938023;
\end{eqnarray*}%
\newline
(ii) $\ $for $h,\ell \geq 1$,
\begin{eqnarray*}
W_{(\alpha ,\beta )}^{(h,\ell )} &\in &\mathfrak{H}_{1}\text{ and }%
W_{(\alpha ,\beta )}\notin \mathfrak{H}_{\infty } \\
&\Leftrightarrow &\sqrt{\tfrac{2}{3}}<\kappa \leq \sqrt{G(h)}=\sqrt{\frac{%
2h^{3}+7h^{2}+8h+3}{2h^{3}+7h^{2}+10h+4}};
\end{eqnarray*}%
\newline
(iii) $\ $for $h,\ell \geq 1$,
\begin{eqnarray*}
W_{(\alpha ,\beta )}^{(h,\ell )} &\in &\mathfrak{H}_{2}\text{ and }%
W_{(\alpha ,\beta )}\notin \mathfrak{H}_{\infty } \\
&\Leftrightarrow &\sqrt{\tfrac{2}{3}}<\kappa \leq \sqrt{\frac{2}{2+b(2,h)}}=%
\sqrt{\frac{8h^{4}+24h^{3}+26h^{2}+12h+2}{8h^{4}+36h^{3}+39h^{2}+18h+3}}.
\end{eqnarray*}
\end{example}

%%%%%%%%%%%%%%%%%%%%%%%%

\section{\label{hypoinv}Hyponormal Invariance Under Powers in the Class $%
\mathcal{A}$}

In this section we study a large class $\mathcal{C}$ of nontrivial pairs of
commuting subnormals such that $W_{(\alpha ,\beta )}\in \mathfrak{H}%
_{1}\bigcap \mathcal{C}\Rightarrow W_{(\alpha ,\beta )}^{(h,\ell )}\in
\mathfrak{H}_{1}\mathcal{\;}$(all $h,\ell \geq 1$). \ The class $\mathcal{C}$
is a subclass of the class $\mathcal{A}$, and it consists of $2$-variable
weighted shifts whose weight diagrams are given in Figure \ref{k-hypo}(i). \
Motivated by the necessary condition for LPCS found in \cite{CuYo2} (see
Lemma \ref{necessary}), we observe that the Berger measure $\xi _{x}$ of the
unilateral weighted shift $W_{x}\equiv \operatorname{shift}\;(\alpha _{00},\alpha
_{10},\cdots )$ admits a unique decomposition as
\begin{equation*}
\xi _{x}\equiv p\delta _{0}+q\delta _{1}+(1-p-q)\rho ,
\end{equation*}%
where $0<p,q<1$, $p+q\leq 1$, and $\rho $ a probability measure with $\rho
(\{0,1\})=0$. \ As a result, a $2$-variable weighted shift $W_{(\alpha
,\beta )}\in \mathcal{C}$ can be parameterized as $W_{(\alpha ,\beta
)}\equiv \left\langle p,q,\rho ,y,a\right\rangle $, with $0<a,y\leq 1$. \ In
Theorem \ref{Thm-main} below we characterize the shifts $W_{(\alpha ,\beta
)} $ which remain hyponormal, $2$-hyponormal or subnormal under the action $%
(h,\ell )\mapsto W_{(\alpha ,\beta )}^{(h,\ell )}\;\;(h,\ell \geq 1)$. \

\begin{theorem}
\label{Thm-main} Let $W_{(\alpha ,\beta )}\equiv \left\langle p,q,\rho
,y,a\right\rangle \in \mathcal{C}$ be the $2$-variable weighted shift whose
weight diagram is given in Figure \ref{k-hypo}(i). \ The following
assertions hold.\newline
(i) \ Assume that $W_{(\alpha ,\beta )}^{(h,\ell )}\in \mathfrak{H}_{1}$
(all $h,\ell \geq 1$). \ Then $0<y\leq m_{1}(a,q):=\sqrt{\frac{q(1-q)}{%
(a^{2}-q)^{2}+q(1-q)}}$.\newline
(ii) \ Assume that $W_{(\alpha ,\beta )}^{(h,\ell )}\in \mathfrak{H}_{2}$ $\
$(all $h,\ell \geq 1$). \ Then $y\leq m_{2}(a,q):=\min \left\{ \sqrt{\frac{%
1-q}{1-a^{2}}},\sqrt{\frac{q}{a^{2}}}\right\} $.\newline
(iii) $\ W_{(\alpha ,\beta )}\in \mathfrak{H}_{\infty }\Longleftrightarrow
y\leq m_{\infty }(a,p,q):=\min \left\{ \sqrt{\frac{p}{1-a^{2}}},\sqrt{\frac{q%
}{a^{2}}}\right\} $.
\end{theorem}

We need an auxiliary lemma, of independent interest.

\begin{lemma}
\label{aux}Let $W_{x}$ be a subnormal unilateral weighted shift, with Berger
measure $\xi _{x}\equiv p\delta _{0}+q\delta _{1}+[1-(p+q)]\rho $, and
recall that $\gamma _{n}$ is the $n$-th moment of $\xi _{x}$, that is, $%
\gamma _{n}=\int s^{n}d\mathcal{\xi }_{x}(s)\;(n\geq 0)$. \ Then $%
\lim_{n\rightarrow \infty }\gamma _{n}=q$.

\begin{proof}
For $n\geq 0$, let $f_{n}(s):=s^{n}\;(0\leq s\leq 1)$. \ Consider the
sequence of nonnegative functions $\{f_{n}\}_{n\geq 0}$. \ Clearly $%
f_{n}\searrow \chi _{\{1\}}$ pointwise, and $\left\vert f_{n}\right\vert
\leq 1\;$(all $n\geq 0$). \ By the Lebesgue Dominated Convergence Theorem,
we have
\begin{equation*}
\lim_{n\rightarrow \infty }\int f_{n}(s)d\rho (s)=\int \chi _{\{1\}}d\rho
(s)=\rho (\{1\})=0
\end{equation*}%
(recall that $\rho (\{0\}\cup \{1\})=0$). \ Thus
\begin{equation*}
\lim_{n\rightarrow \infty }\gamma _{n}=\lim_{n\rightarrow \infty }\int s^{n}d%
\mathcal{\xi }_{x}(s)=\lim_{n\rightarrow \infty }\int f_{n}(s)d\mathcal{\xi }%
_{x}(s)=q+[1-(p+q)]\cdot \lim_{n\rightarrow \infty }\int f_{n}(s)d\rho (s)=q,
\end{equation*}%
as desired.
\end{proof}
\end{lemma}

\begin{proof}[Proof of Theorem \protect\ref{Thm-main}]
For fixed $h,\ell \geq 1$, $0\leq m\leq h-1$ and $0\leq n\leq \ell -1$, we
recall that $\mathcal{H}_{(m,n)}=\bigvee_{i,j=0}^{\infty }\{e_{(hi+m,\ell
j+n)}:h,\ell \geq 1\}$ and $\ell ^{2}(\mathbb{Z}_{+}^{2})\equiv
\bigoplus_{m=0}^{h-1}\bigoplus_{n=0}^{\ell -1}\mathcal{H}_{(m,n)}$. \ For $%
W_{(\alpha ,\beta )}^{(h,\ell )}|_{\mathcal{H}_{(0,0)}}$, we refer to the
weight diagram in Figure \ref{k-hypo}(ii). \ In the decomposition $\mathcal{%
\xi }_{x}\equiv p\delta _{0}+q\delta _{1}+[1-(p+q)]\rho $, we may assume,
without loss of generality, that $q<1$; for, the condition $q=1$ and
hyponormality immediately imply the subnormality of $W_{(\alpha ,\beta )}$.

Given $h,\ell \geq 1$, we consider the moments associated with $W_{(\alpha
,\beta )}^{(h,\ell )}$ of order $\mathbf{k}$, that is,%
\begin{equation}
\gamma _{\mathbf{k}}(W_{(\alpha ,\beta )}^{(h,\ell )})=\left\{
\begin{tabular}{ll}
$1,$ & $\text{if }k_{1}=0\text{ and }k_{2}=0$ \\
$\gamma _{k_{1}h}(W_{(\alpha ,\beta )}),$ & $\text{if }k_{1}\geq 1\text{ and
}k_{2}=0$ \\
$y^{2},$ & $\text{if }k_{1}=0\text{ and }k_{2}\geq 1$ \\
$a^{2}y^{2},$ & $\text{if }k_{1}\geq 1\text{ and }k_{2}\geq 1$.%
\end{tabular}%
\right.  \label{moment}
\end{equation}%
(i) \ From Lemma \ref{thmbackext}, we observe that for $h,\ell \geq 1$, $%
W_{(\alpha ,\beta )}^{(h,\ell )}\mathbf{|}_{\mathcal{M}_{1}}\cong (I\otimes
S_{a},U_{+}\otimes I)\in \mathfrak{H}_{\infty }$ and $W_{(\alpha ,\beta
)}^{(h,\ell )}\mathbf{|}_{\mathcal{N}_{1}}\in \mathfrak{H}_{\infty }$. \
Thus, by Lemma \ref{khypo}, to verify the hyponormality of $W_{(\alpha
,\beta )}^{(h,\ell )}$, it suffices to apply the Six-point Test (Lemma \ref%
{joint hypo}) to $W_{(\alpha ,\beta )}^{(h,\ell )}$ at $\mathbf{k}=(0,0)$. \
We then have%
\begin{eqnarray*}
M_{(0,0)}(1)(W_{(\alpha ,\beta )}^{(h,\ell )}) &=&\left(
\begin{array}{ccc}
1 & \gamma _{h}(W_{(\alpha ,\beta )}) & y^{2} \\
\gamma _{h}(W_{(\alpha ,\beta )}) & \gamma _{2h}(W_{(\alpha ,\beta )}) &
a^{2}y^{2} \\
y^{2} & a^{2}y^{2} & y^{2}%
\end{array}%
\right) \geq 0\;\;\text{(all }h\geq 1\text{)} \\
&\Rightarrow &H:=\left(
\begin{array}{ccc}
1 & q & y^{2} \\
q & q & a^{2}y^{2} \\
y^{2} & a^{2}y^{2} & y^{2}%
\end{array}%
\right) \geq 0\;\;\text{(by Lemma \ref{aux})} \\
&\Longleftrightarrow &0<y\leq m_{1}(a,q):=\sqrt{\frac{q(1-q)}{%
(a^{2}-q)^{2}+q(1-q)}},
\end{eqnarray*}%
as desired. \ Observe that the function $m_{1}$ satisfies the following
properties:\newline
(i$_{1}$) \ $0<m_{1}(a,q)\leq 1$ on the square $\left( 0,1\right] \times
(0,1)$; \newline
(i$_{2}$) $\lim_{q\rightarrow 0^{+}}\;m_{1}(a,q)=0=\lim_{q\rightarrow
1^{+}}\;m_{1}(a,q)=0\;\;$(for all $a$); \newline
(i$_{3}$) $\lim_{a\rightarrow 0^{+}}\;m_{1}(a,q)=\sqrt{1-q}\;\;$(for all $q$%
); \newline
(i$_{4}$) $m_{1}(1,q)=\sqrt{q}\;\;$(for all $q$); and \newline
(i$_{5}$) $m_{1}(a,q)=1\Longleftrightarrow q=a^{2}$. \ \newline
Thus near the edges of the square the hyponormality of $W_{(\alpha ,\beta
)}^{(h,\ell )}$ for all $h$ and $\ell $ forces $y$ to be small, while along
the parabola $q=a^{2}$ the values of $y$ can reach $1$.

(ii) \ From Lemmas \ref{thmbackext} and \ref{khypo}, and the fact that $%
W_{(\alpha ,\beta )}^{(h,\ell )}\mathbf{|}_{\mathcal{M}_{1}}$, $W_{(\alpha
,\beta )}^{(h,\ell )}\mathbf{|}_{\mathcal{N}_{1}}\in \mathfrak{H}_{\infty }$%
, to verify the $2$-hyponormality of $W_{(\alpha ,\beta )}^{(h,\ell )}$ ($%
h,\ell \geq 1$) it suffices to apply the $15$-point Test to $W_{(\alpha
,\beta )}^{(h,\ell )}$ at $\mathbf{k}=(0,0)$. \ By direct computation (i.e.,
interchanging rows and columns, and discarding some redundant rows and
columns), it is straightforward to observe that the positivity of the $%
10\times 10$ matrix $M_{(0,0)}(2)(W_{(\alpha ,\beta )}^{(h,\ell )})$ is
determined by that of the following $5\times 5$ matrix:
\begin{equation*}
P(h):=\left(
\begin{array}{cc}
\left(
\begin{array}{cc}
1 & \gamma _{h}(W_{(\alpha ,\beta )}) \\
\gamma _{h}(W_{(\alpha ,\beta )}) & \gamma _{2h}(W_{(\alpha ,\beta )})%
\end{array}%
\right) &
\begin{pmatrix}
y^{2} & \gamma _{2h}(W_{(\alpha ,\beta )}) & a^{2}y^{2} \\
a^{2}y^{2} & \gamma _{3h}(W_{(\alpha ,\beta )}) & a^{2}y^{2}%
\end{pmatrix}
\\
\begin{pmatrix}
y^{2} & a^{2}y^{2} \\
\gamma _{2h}(W_{(\alpha ,\beta )}) & \gamma _{3h}(W_{(\alpha ,\beta )}) \\
a^{2}y^{2} & a^{2}y^{2}%
\end{pmatrix}
&
\begin{pmatrix}
y^{2} & a^{2}y^{2} & a^{2}y^{2} \\
a^{2}y^{2} & \gamma _{4h}(W_{(\alpha ,\beta )}) & a^{2}y^{2} \\
a^{2}y^{2} & a^{2}y^{2} & a^{2}y^{2}%
\end{pmatrix}%
\end{array}%
\right) .
\end{equation*}%
Thus the assumption $W_{(\alpha ,\beta )}^{(h,\ell )}\in \mathfrak{H}_{2}$ $%
\ $(all $h,\ell \geq 1$) readily implies that
\begin{equation*}
P\equiv P(\infty ):=\left(
\begin{array}{cc}
\left(
\begin{array}{cc}
1 & q \\
q & q%
\end{array}%
\right) &
\begin{pmatrix}
y^{2} & q & a^{2}y^{2} \\
a^{2}y^{2} & q & a^{2}y^{2}%
\end{pmatrix}
\\
\begin{pmatrix}
y^{2} & a^{2}y^{2} \\
q & q \\
a^{2}y^{2} & a^{2}y^{2}%
\end{pmatrix}
&
\begin{pmatrix}
y^{2} & a^{2}y^{2} & a^{2}y^{2} \\
a^{2}y^{2} & q & a^{2}y^{2} \\
a^{2}y^{2} & a^{2}y^{2} & a^{2}y^{2}%
\end{pmatrix}%
\end{array}%
\right) \geq 0\text{ (using Lemma \ref{aux}).}
\end{equation*}%
Since$\left(
\begin{array}{cc}
1 & q \\
q & q%
\end{array}%
\right) $ is positive and invertible, we can apply Smul'jan Lemma (Lemma \ref%
{smu}) to $P$:
\begin{eqnarray*}
P &\geq &0\Longleftrightarrow
\begin{pmatrix}
\frac{y^{2}\left( q\left( 1-(1-2a^{2})y^{2}\right) -q^{2}-a^{4}y^{2}\right)
}{(1-q)q} & 0 & a^{2}y^{2}-\frac{a^{4}y^{4}}{q} \\
0 & 0 & 0 \\
a^{2}y^{2}-\frac{a^{4}y^{4}}{q} & 0 & a^{2}y^{2}-\frac{a^{4}y^{4}}{q}%
\end{pmatrix}%
\geq 0 \\
&\Longleftrightarrow &\left\{
\begin{array}{ccc}
(1-a^{2})y^{2} & \leq & 1-q \\
a^{2}y^{2} & \leq & q%
\end{array}%
\right. \Leftrightarrow y\leq \min \left\{ \sqrt{\frac{1-q}{1-a^{2}}},\sqrt{%
\frac{q}{a^{2}}}\right\} \text{.}
\end{eqnarray*}%
Therefore, $W_{(\alpha ,\beta )}^{(h,\ell )}\in \mathfrak{H}_{2}\;\;$(all $%
h,\ell $)$\;\Rightarrow y\leq \min \left\{ \sqrt{\frac{1-q}{1-a^{2}}},\sqrt{%
\frac{q}{a^{2}}}\right\} $, as desired. \ (The reader will notice that $\min
\left\{ \sqrt{\frac{1-q}{1-a^{2}}},\sqrt{\frac{q}{a^{2}}}\right\} \leq 1$;
for, if $q>a^{2}$ then $1-q<1-a^{2}$.)\newline
(iii) \ From Figure \ref{k-hypo}(i), we observe that $W_{(\alpha ,\beta )}%
\mathbf{|}_{\mathcal{M}}$ is subnormal with Berger measure%
\begin{equation*}
\mu _{\mathcal{M}}=\left( (1-a^{2})\delta _{0}+a^{2}\delta _{1}\right)
\times \delta _{1}.
\end{equation*}%
Note that $\left\Vert \frac{1}{t}\right\Vert _{L^{1}(\mu _{\mathcal{M}})}=1$
and $(\mu _{\mathcal{M}})_{ext}^{X}=(1-a^{2})\delta _{0}+a^{2}\delta _{1}$.
\ We now apply Lemma \ref{backext} to the $2$-variable weighted shift $%
W_{(\alpha ,\beta )}$ and to the subspace $\mathcal{M}$. \ It follows that
the necessary and sufficient condition for $W_{(\alpha ,\beta )}$ to be
subnormal is
\begin{equation*}
y^{2}\left\Vert \frac{1}{t}\right\Vert _{L^{1}(\mu _{\mathcal{M}})}(\mu _{%
\mathcal{M}})_{ext}^{X}\leq \mathcal{\xi }_{x}=p\delta _{0}+q\delta
_{1}+[1-(p+q)]\rho ,
\end{equation*}%
or equivalently,%
\begin{equation*}
\left\{
\begin{array}{ccc}
(1-a^{2})y^{2} & \leq & p \\
a^{2}y^{2} & \leq & q.%
\end{array}%
\right.
\end{equation*}%
Thus we have the desired result. \ The proof of the theorem is now complete.
\end{proof}

\begin{remark}
\label{Re4} (i) \ By a direct calculation, we note that
\begin{equation*}
\begin{tabular}{l}
$a^{2}\leq q\Leftrightarrow \frac{1-q}{1-a^{2}}\leq \frac{q(1-q)}{%
(a^{2}-q)^{2}+q(1-q)}$ and $a^{2}>q\Leftrightarrow \frac{q}{a^{2}}<\frac{%
q(1-q)}{(a^{2}-q)^{2}+q(1-q)}$.%
\end{tabular}%
\end{equation*}%
Thus it is always true that $\min \left\{ \sqrt{\frac{q}{a^{2}}},\sqrt{\frac{%
1-q}{1-a^{2}}}\right\} \leq \sqrt{\frac{q(1-q)}{(a^{2}-q)^{2}+q(1-q)}}$.%
\newline
(ii) \ For $h\geq 2$ and $k,\ell \geq 1$, the necessary conditions in (i)
and (ii) in Theorem \ref{Thm-main} are not sufficient for power invariance.
\ To show this, we let $W_{(\alpha ,\beta )}$ denote the $2$-variable
weighted shift whose weight diagram is given in Figure \ref{k-hypo}(i), with
$a=\sqrt{\frac{1}{2}}$ and $W_{x}$ is as in (\ref{alpha}). \ Furthermore,
for given small $\varepsilon >0$ and $h\geq 1$, we let
\begin{equation}
\kappa :=\left\{
\begin{tabular}{ll}
$\sqrt{\frac{20}{23}}$, & if $k=1$ \\
$\sqrt{\frac{2}{2+b(k,h)}}+\varepsilon (h)$, & if $k\geq 2$%
\end{tabular}%
\right\}  \label{con0}
\end{equation}%
provided that $\sqrt{\frac{2}{2+b(k,h)}}+\varepsilon (h)<1$. \ (Since, for $%
k,h\geq 1$, $b(k,h)>0$, $\lim_{k\rightarrow \infty }b(k,h)=1$ and $%
\lim_{h\rightarrow \infty }b(k,h)=0$, it is possible to choose $\kappa $
given in (\ref{con0})). \ By Theorem \ref{Thm(invariant)-00}, we note that
for $h\geq 2$, $k=1,2$ and $\ell \geq 1$, $W_{(\alpha ,\beta )}^{(h,\ell
)}\notin \mathfrak{H}_{k}$ ($k=1,2$). \ Observe that
\begin{equation}
\begin{tabular}{l}
$\left\{
\begin{tabular}{ll}
$0<\kappa \leq \sqrt{\frac{q(1-q)}{a^{4}+q-2a^{2}q}}\Leftrightarrow \kappa
\leq 1$, & if $k=1$ \\
&  \\
$\left\{ a^{2}\kappa ^{2}\leq q\leq \left( 1-\kappa ^{2}\right) +a^{2}\kappa
^{2}\right\} \Leftrightarrow \kappa \leq 1$, & if $k\geq 2.$%
\end{tabular}%
\right. $%
\end{tabular}
\label{con}
\end{equation}%
If we choose $\kappa $ given in (\ref{con0}), then (\ref{con}) is always
true. \ Thus the $2$-variable weighted shift $W_{(\alpha ,\beta )}$ given in
Figure \ref{k-hypo}(i) satisfies the necessary conditions in Theorem \ref%
{Thm-main}, but for $h\geq 2$ and $\ell \geq 1$, $W_{(\alpha ,\beta
)}^{(h,\ell )}\notin \mathfrak{H}_{1}$.
\end{remark}

Throughout this section, we have focused on the question of hyponormality
for shifts in the class $\mathcal{A}$. \ We now turn our attention to $2$%
-hyponormality, in the hope of detecting to what extent one can expect
invariance under powers in this class. \ Along the way we will discover that
there is a large subclass, $\mathcal{S}_{1}$, for which things work
extremely well. \ We will show this in Section \ref{S1inv}. 

As we saw in Section \ref{Sect4}, for a general operator $T$ on Hilbert
space and for all $m\geq 2$, we know that the $k$-hyponormality of $T^{m} \ %
(k\geq 2)$ need not imply the $k$-hyponormality of $T$. \ But it is still unknown
whether the $k$-hyponormality of $T$ $(k\geq 2)$ implies the $k$%
-hyponormality of $T^{m}$ $(m\geq 2)$, even when $T$ is a weighted shift
(see Problem \ref{question0}). \ We now show that there exists a subclass of $2$%
-variable weighted shifts $W_{(\alpha ,\beta )}\in \mathcal{A}$ for which the
$2$-hyponormality of $W_{(\alpha ,\beta )}$ does imply the $2$%
-hyponormality of $W_{(\alpha ,\beta )}^{(2,1)}$.

The motivation behind the construction in Theorem \ref{Thm8} comes from
Figure \ref{range}. \ Indeed, inspection of the values of $a$ that
illustrate the gap between the hyponormality of $W_{(\alpha ,\beta )}$ and
that of its powers suggests that something similar may work for $k$%
-hyponormality. \ We saw in Theorem \ref{Thm(invariant)-00} that a value
smaller than $a_{int}$ (namely, $a=\sqrt{\frac{1}{2}}$) did the job in
separating the $k$-hyponormality of the powers from the $k$-hyponormality of
the pair. \ For the subclass below, however, we do see a propagation effect, that is, the $(2,1)$ power is $2$-hyponormal whenever the 
original weighted shift is.

\begin{theorem}
\label{Thm8}Let $W_{(\alpha ,\beta )}\equiv (T_{1},T_{2})$ be the $2$%
-variable weighted shift whose weight diagram is given in Figure \ref{k-hypo}%
(i) (where $0<a \leq \sqrt {\frac{1}{2}}$ and $W_{x}$ is as in (\ref{alpha})). \ Then%
\newline
(i) $\ W_{(\alpha ,\beta )}\in \mathfrak{H}_{0}$;\newline
(ii) \ $W_{(\alpha ,\beta )}\in \mathfrak{H}_{2}\Leftrightarrow \kappa
\leq h_{2}(a) :=\sqrt{\frac {9(9-16a^2)}{157-360a^2+144a^4}}$;\newline
(iii) \ $W_{(\alpha ,\beta )}^{(2,1)}\in \mathfrak{H}_{2}\Leftrightarrow \kappa \leq h_{2}^{(2,1)}(a):=\sqrt{\frac{225(15-28a^2)}
{6238-15015a^2+6300a^4}}$;\newline
(iv) \ $h_{2}(a)<h_{2}^{(2,1)}(a)$ for $a \in (0,\sqrt{\frac{1}{2}}]$.
\end{theorem}

\begin{proof}
(i) \ This follows easily once we know that $W_{x}$ is subnormal.\newline
(ii) \ This is part (iii) of Theorem \ref{examplethm}.\newline
(iii) \ Recall that for $n=0,1$, $\mathcal{H}_{n}\equiv
\bigvee_{i=0}^{\infty }\{e_{(2i+n,j)}:j=0,1,2,\cdots \}$ and $\ell ^{2}(%
\mathbb{Z}_{+}^{2})\equiv \mathcal{H}_{0}\bigoplus \mathcal{H}_{1}$. \ Note
that%
\begin{equation*}
W_{(\alpha ,\beta )}^{(2,1)}\equiv (T_{1}^{2},T_{2})\cong W_{(\alpha ,\beta
)}^{(2,1)}|_{{\mathcal{H}_{0}}}\bigoplus W_{(\alpha ,\beta )}^{(2,1)}|_{{%
\mathcal{H}_{1}}}.
\end{equation*}%
By Lemma \ref{thmbackext}, we observe that $W_{(\alpha ,\beta )}^{(2,1)}|_{{%
\mathcal{H}_{1}}}\cong (W_{\alpha (2:1)}\oplus (I\otimes U_{+}),T_{2}|_{{%
\mathcal{H}_{1}}})$ is subnormal, because $W_{\alpha (2:1)}=\operatorname{shift}%
\;(x_{1}x_{2},x_{3}x_{4},\cdots )$ has an atom at $\{1\}$. \ Hence $%
W_{(\alpha ,\beta )}^{(2,1)}\in \mathfrak{H}_{2}$ if and only if $W_{(\alpha
,\beta )}^{(2,1)}|_{{\mathcal{H}_{0}}}:=\mathbf{(}T_{1}^{2},T_{2})|_{%
\mathcal{H}_{0}}\in \mathfrak{H}_{2}$. \ Let $\mathcal{M}_{1}(0)$ (resp. $%
\mathcal{N}_{1}(0)$) be the subspace of $\mathcal{H}_{0}$ spanned by
canonical orthonormal basis vectors associated to indices $\mathbf{k}%
=(k_{1},k_{2})$ with $k_{1}\geq 0$ and $k_{2}\geq 1$ (resp. $k_{1}\geq 1$
and $k_{2}\geq 0$). \ By Lemma \ref{thmbackext}, we note that $W_{(\alpha
,\beta )}^{(2,1)}|_{\mathcal{H}_{0}}$ on $\mathcal{M}_{1}(0)$ (resp. $%
\mathcal{N}_{1}(0)$) is subnormal, because $\operatorname{shift}%
\;(x_{2}x_{3},x_{4}x_{5},\cdots )\ $(resp. $\operatorname{shift}\;(y,1,1\cdots )$)
has an atom at $\{1\}$. $\ $Thus, by Lemma \ref{khypo}, to verify the $2$%
-hyponormality of $W_{(\alpha ,\beta )}^{(2,1)}$ it suffices to apply the $%
15 $-point Test to $W_{(\alpha ,\beta )}^{(2,1)}|_{\mathcal{H}_{0}}$ at $%
\mathbf{k}=(0,0)$. \ Note that the moments associated with $\mathbf{(}%
T_{1}^{2},T_{2})|_{\mathcal{H}_{0}}$ of order $\mathbf{k}$ are
\begin{equation*}
\gamma _{\mathbf{k}}(W_{(\alpha ,\beta )}^{(2,1)}|_{\mathcal{H}%
_{0}})=\left\{
\begin{tabular}{ll}
$1,$ & $\text{if }k_{1}=0\text{ and }k_{2}=0$ \\
$\prod\limits_{i=1}^{k_{1}}x_{2(i-1)}^{2}x_{2i-1}^{2},$ & $\text{if }%
k_{1}\geq 1\text{ and }k_{2}=0$ \\
$\kappa ^{2}$ & $\text{if }k_{1}=0\text{ and }k_{2}\geq 1$ \\
$a^2 \kappa ^{2},$ & $\text{if }k_{1}\geq 1\text{ and }k_{2}\geq
1 $%
\end{tabular}%
\right. .
\end{equation*}%
Since the third and sixth rows of $M_{(0,0)}(2)|_{\mathcal{H}_{0}}$ are
identical, if we multiply $\frac{1}{y^{2}}$ and then apply row and column
operations to $M_{(0,0)}(2)|_{\mathcal{H}_{0}}$, then we have $%
M_{(0,0)}(2)|_{\mathcal{H}_{0}}\geq 0\Leftrightarrow \tilde{M}_{(0,0)}\geq 0$%
, where
\begin{equation*}
\tilde{M}_{(0,0)}:=\left(
\begin{array}{cc}
\begin{pmatrix}
\frac{3}{5} & \frac{4}{7} & a^{2} \\
\frac{4}{7} & \frac{5}{9} & a^{2} \\
a^{2} & a^{2} & 1%
\end{pmatrix}
&
\begin{pmatrix}
\frac{2}{3} & a^{2} \\
\frac{3}{5} & a^{2} \\
1 & a^{2}%
\end{pmatrix}
\\
\begin{pmatrix}
\frac{2}{3} & \frac{3}{5} & 1 \\
a^{2} & a^{2} & a^{2}%
\end{pmatrix}
&
\begin{pmatrix}
\frac{1}{y^{2}} & a^{2} \\
a^{2} & a^{2}%
\end{pmatrix}%
\end{array}%
\right) .
\end{equation*}

We now apply the Nested Determinant Test to $\tilde{M}_{(0,0)}$; let $d_{p}$ denote the determinant of the $p$ by $p$ principal minor consisting of the first $p$ rows and columns of $\tilde{M}_{(0,0)}$. \ Thus, $d_{1}=\frac{1}{\kappa^2} > 0$ and $d_{2}=\frac{27-20\kappa^2}{45\kappa^2}$. \ Since $\kappa <1$, it follows easily that $d_{2} >0$. \ Now, $d_{3} = \frac{27 - 45 a^4 - 47 \kappa^2 + 60 a^2 \kappa^2}{45 \kappa^2}$, and we conclude that $d_{3}>0$ if $\kappa <\kappa_{3}:=3 \sqrt{\frac{3 - 5 a^4}{47 - 60 a^2}}$. \ Similarly, $d_{4} =\frac{3375 - 6300 a^4 - 6238 \kappa^2
 + 6510 a^2 \kappa^2 + 2205 a^4 \kappa^2}{496125 \kappa^2}$ and $d_{4}>0$ whenever $\kappa < \kappa_{4}:= 15 \sqrt{\frac{15 - 28 a^4}
 {6238 - 6510 a^2 - 2205 a^4}}$. \  Finally, 
 \begin{equation*}
 d_{5}=\frac{3375 a^2 - 9675 a^4 + 6300 a^6 - 6238 a^2 \kappa^2 + 
 21253 a^4 \kappa^2 - 21315 a^6 \kappa^2 + 
 6300 a^8 \kappa^2}{496125 \kappa^2},
 \end{equation*}
 and in this case we can guarantee the nonnegativity of $d_{5}$ if $\kappa \leq \kappa_{5}:=
 15 \sqrt{\frac{15 - 43 a^2 + 28 a^4}{6238 - 21253 a^2 + 21315 a^4 - 
 6300 a^6}}$. \ Visual inspection of the graphs of $\kappa_{3}$, $\kappa_{4}$, $\kappa_{5}$ and $h_{2}$ on the common interval $[0,\sqrt{\frac{1}{2}}]$ reveals that $h_{2}<\kappa_{5}<\kappa_{4}<\kappa_{3}$. \ We thus conclude that $W_{(\alpha ,\beta )}^{(2,1)}$ is $2$-hyponormal if and only if $\kappa \leq \kappa_{5} \equiv h_{2}^{(2,1)}$, as desired. 
\end{proof}

\begin{corollary}
\label{re2} Let $W_{(\alpha ,\beta)}$ be as in Theorem \ref{Thm8}, let $a \in (0,\sqrt{\frac{1}{2}}]$, and assume that 
$W_{(\alpha ,\beta
)}\in \mathfrak{H}_{2}$. \ Then $W_{(\alpha ,\beta )}^{(2,1)}\in \mathfrak{H}%
_{2}$.
\end{corollary}

\begin{proof}
This is straightforward from Theorem \ref{Thm8}(iv).
\end{proof}

\begin{remark}
Looking at Theorem \ref{Thm8}, if seems natural to conjecture that a similar
result should work for $k$-hyponormality ($k>2$). \ That is, perhaps one has 
$W_{(\alpha ,\beta )}\in \mathfrak{H}_{k}\Longrightarrow W_{(\alpha ,\beta )}^{(2,1)}$ whenever $0<a\leq \sqrt{\frac{1}{2}}$. 
\end{remark}

%%%%%%%%%%%%%%%%%%%%%%%%%%%

\section{\label{S1inv}The Class $\mathcal{S}_{1}$ Is Invariant Under All
Powers}

In Section \ref{hypoinv} we dealt with $2$-variable weighted shifts of the
form $W_{(\alpha ,\beta )}\equiv \left\langle p,q,\rho ,y,a\right\rangle $
and established some results about hyponormality, $2$-hyponormality and
subnormality. \ We now restrict attention to the case $\rho =0$, and assume
that $W_{(\alpha ,\beta )}\in \mathfrak{H}_{1}$; that is, $W_{(\alpha ,\beta
)}\in \mathcal{S}_{1}$. \ Under this assumption, we will now sharpen the
hyponormality results. \ Recall that, without loss of generality, every $2$%
-variable weighted shift $W_{(\alpha ,\beta )}\in \mathcal{S}_{1}$ is
completely determined by the three parameters $x:=\alpha _{(0,0)}$, $%
y:=\beta _{(0,0)}$ and $a:=\alpha _{(0,1)}$; cf. Figure \ref{S1}(i). \ As
before, we shall denote such a shift by $\left\langle x,y,a\right\rangle $;
of course, we always assume $0<x,y,a\leq 1$, and moreover $ay\leq x$ (since
we need to ensure that $\operatorname{shift}\;(\beta _{10},\beta _{11},\beta
_{12},\cdots )\equiv \operatorname{shift}\;(\frac{ay}{x},1,1,\cdots )$ is subnormal).

\setlength{\unitlength}{1mm} \psset{unit=1mm}
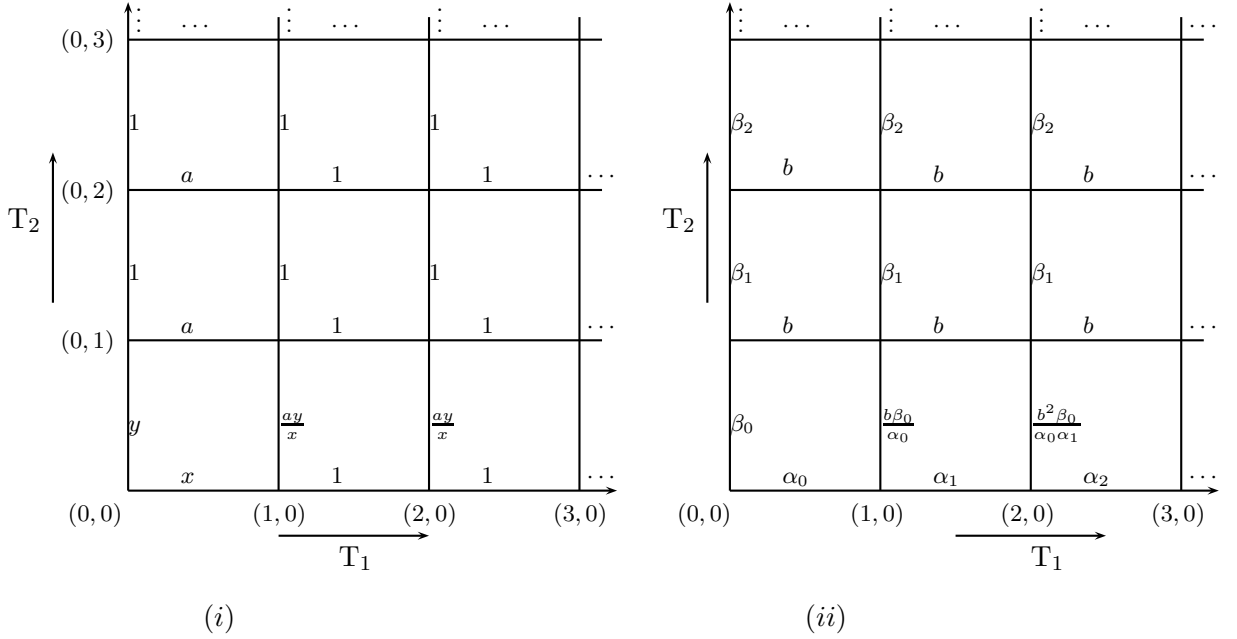
\begin{figure}[th]
\begin{center}
\begin{picture}(165,85)

\psline{->}(20,20)(85,20)
\psline(20,40)(83,40)
\psline(20,60)(83,60)
\psline(20,80)(83,80)
\psline{->}(20,20)(20,85)
\psline(40,20)(40,83)
\psline(60,20)(60,83)
\psline(80,20)(80,83)

\put(12,16){\footnotesize{$(0,0)$}}
\put(36.5,16){\footnotesize{$(1,0)$}}
\put(56.5,16){\footnotesize{$(2,0)$}}
\put(76.5,16){\footnotesize{$(3,0)$}}

\put(27,21){\footnotesize{$x$}}
\put(47,21){\footnotesize{$1$}}
\put(67,21){\footnotesize{$1$}}
\put(81,21){\footnotesize{$\cdots$}}

\put(27,41){\footnotesize{$a$}}
\put(47,41){\footnotesize{$1$}}
\put(67,41){\footnotesize{$1$}}
\put(81,41){\footnotesize{$\cdots$}}

\put(27,61){\footnotesize{$a$}}
\put(47,61){\footnotesize{$1$}}
\put(67,61){\footnotesize{$1$}}
\put(81,61){\footnotesize{$\cdots$}}

\put(27,81){\footnotesize{$\cdots$}}
\put(47,81){\footnotesize{$\cdots$}}
\put(67,81){\footnotesize{$\cdots$}}

\psline{->}(40,14)(60,14)
\put(48,10){$\rm{T}_1$}
\psline{->}(10,45)(10,65)
\put(4,55){$\rm{T}_2$}

\put(11,39){\footnotesize{$(0,1)$}}
\put(11,59){\footnotesize{$(0,2)$}}
\put(11,79){\footnotesize{$(0,3)$}}

\put(20,28){\footnotesize{$y$}}
\put(20,48){\footnotesize{$1$}}
\put(20,68){\footnotesize{$1$}}
\put(21,81){\footnotesize{$\vdots$}}

\put(40,28){\footnotesize{$\frac{ay}{x}$}}
\put(40,48){\footnotesize{$1$}}
\put(40,68){\footnotesize{$1$}}
\put(41,81){\footnotesize{$\vdots$}}

\put(60,28){\footnotesize{$\frac{ay}{x}$}}
\put(60,48){\footnotesize{$1$}}
\put(60,68){\footnotesize{$1$}}
\put(61,81){\footnotesize{$\vdots$}}

\put(30,2){$(i)$}

%\\\\\\\\\

\put(110,2){$(ii)$}

\psline{->}(130,14)(150,14) \put(140,10){$\rm{T}_1$}
\psline{->}(97,45)(97,65) \put(91,55){$\rm{T}_2$}

\psline{->}(100,20)(165,20)
\psline(100,40)(163,40)
\psline(100,60)(163,60)
\psline(100,80)(163,80)

\psline{->}(100,20)(100,85)
\psline(120,20)(120,83)
\psline(140,20)(140,83)
\psline(160,20)(160,83)

\put(93,16){\footnotesize{$(0,0)$}}
\put(116,16){\footnotesize{$(1,0)$}}
\put(136,16){\footnotesize{$(2,0)$}}
\put(156,16){\footnotesize{$(3,0)$}}

\put(107,21){\footnotesize{$\alpha_{0}$}}
\put(127,21){\footnotesize{$\alpha_{1}$}}
\put(147,21){\footnotesize{$\alpha_{2}$}}
\put(161,21){\footnotesize{$\cdots$}}

\put(107,41){\footnotesize{$b$}}
\put(127,41){\footnotesize{$b$}}
\put(147,41){\footnotesize{$b$}}
\put(161,41){\footnotesize{$\cdots$}}

\put(107,62){\footnotesize{$b$}}
\put(127,61){\footnotesize{$b$}}
\put(147,61){\footnotesize{$b$}}
\put(161,61){\footnotesize{$\cdots$}}

\put(107,81){\footnotesize{$\cdots$}}
\put(127,81){\footnotesize{$\cdots$}}
\put(147,81){\footnotesize{$\cdots$}}
\put(161,81){\footnotesize{$\cdots$}}

\put(100,28){\footnotesize{$\beta_{0}$}}
\put(100,48){\footnotesize{$\beta_{1}$}}
\put(100,68){\footnotesize{$\beta_{2}$}}
\put(101,81){\footnotesize{$\vdots$}}

\put(120,28){\footnotesize{$\frac{b\beta_{0}}{\alpha_{0}}$}}
\put(120,48){\footnotesize{$\beta_{1}$}}
\put(120,68){\footnotesize{$\beta_{2}$}}
\put(121,81){\footnotesize{$\vdots$}}

\put(140,28){\footnotesize{$\frac{b^2\beta_{0}}{\alpha_{0}\alpha_{1}}$}}
\put(140,48){\footnotesize{$\beta_{1}$}}
\put(140,68){\footnotesize{$\beta_{2}$}}
\put(141,81){\footnotesize{$\vdots$}}
\end{picture}
\end{center}
\caption{Weight diagram of a generic $2$-variable weighted shift in $%
\mathcal{S}_{1}$ and the $2$-variable weighted shift in Lemma \protect\ref%
{thmbackext}, respectively.}
\label{S1}
\end{figure}

First, we wish to obtain a canonical representation for the powers $%
\left\langle x,y,a\right\rangle ^{(h,\ell )}$ as an orthogonal direct sum of
$2$-variable weighted shifts in $\mathcal{S}_{1}$. \ In what follows, we
abbreviate the orthogonal direct sums of $m$ copies of a shift $\left\langle
x,y,a\right\rangle $ by $m\cdot \left\langle x,y,a\right\rangle $.

\begin{proposition}
\label{representation}Let $\left\langle x,y,a\right\rangle \in \mathcal{S}%
_{1}$ and let $h,\ell \geq 1$. \ Then
\begin{equation*}
\left\langle x,y,a\right\rangle ^{(h,\ell )}\cong \left\langle
x,y,a\right\rangle \bigoplus (h-1)\cdot \left\langle 1,\frac{ay}{x}%
,1\right\rangle \bigoplus (\ell -1)\cdot \left\langle a,1,a\right\rangle
\bigoplus (h-1)(\ell -1)\cdot \left\langle 1,1,1\right\rangle .
\end{equation*}
\end{proposition}

\begin{proof}
We decompose the space $\ell ^{2}(\mathbb{Z}_{+}^{2})$ as the orthogonal
direct sum of $h\ell $ subspaces $\mathcal{H}_{(m,n)}$, each isometrically
isomorphic to $\ell ^{2}(\mathbb{Z}_{+}^{2})$, namely $\mathcal{H}%
_{(m,n)}:=\bigvee_{i,j=0}^{\infty }\{e_{(hi+m,\ell j+n)}\}\;\;(0\leq m\leq
h-1,0\leq n\leq \ell -1)$. \ This particular decomposition allows us to
write the power $\left\langle x,y,a\right\rangle ^{(h,\ell )}$ as the
orthogonal direct sum $\bigoplus_{0\leq m\leq h-1,0\leq n\leq \ell
-1}\left\langle x,y,a\right\rangle ^{(h,\ell )}|_{\mathcal{H}_{(m,n)}}$. \
We will now identify each of the summands $\left\langle x,y,a\right\rangle
^{(h,\ell )}|_{\mathcal{H}_{(m,n)}}\;\;(0\leq m\leq h-1,0\leq n\leq \ell -1)$%
. $\ $

\textbf{Case 1}: \ ($m=0$, $n=0$) \ Direct inspection of the weight families
$\alpha $ and $\beta $ shows that
\begin{equation*}
\left\langle x,y,a\right\rangle ^{(h,\ell )}e_{(hi,\ell j)}=\left\langle
x,y,a\right\rangle e_{(hi,\ell j)},
\end{equation*}%
and therefore
\begin{equation*}
\left\langle x,y,a\right\rangle ^{(h,\ell )}|_{\mathcal{H}_{(0,0)}}\cong
\left\langle x,y,a\right\rangle .
\end{equation*}

\textbf{Case 2}: \ ($m>0$, $n=0$) \ In this case the generic basis vector of
$\mathcal{H}_{(m,0)}$is $e_{(hi+m,\ell j)}$, so that $\left\langle
x,y,a\right\rangle ^{(h,\ell )}e_{(hi+m,\ell j)}=\left\langle 1,\frac{ay}{x}%
,1\right\rangle e_{(hi+m,\ell j)}$. \ It follows that $\left\langle
x,y,a\right\rangle ^{(h,\ell )}|_{\mathcal{H}_{(m,0)}}\cong \left\langle 1,%
\frac{ay}{x},1\right\rangle $.

\textbf{Case 3}: \ ($m=0$, $n>0$) \ In this case the generic basis vector of
$\mathcal{H}_{(0,n)}$is $e_{(hi,\ell j+n)}$, and therefore $\left\langle
x,y,a\right\rangle ^{(h,\ell )}e_{(hi,\ell j+n)}=\left\langle
a,1,a\right\rangle e_{(hi,\ell j+n)}$. \ It follows that $\left\langle
x,y,a\right\rangle ^{(h,\ell )}|_{\mathcal{H}_{(0,n)}}\cong \left\langle
a,1,a\right\rangle $.

\textbf{Case 4}: \ ($m>0$, $n>0$) \ Since $\mathcal{H}_{(m,n)}\subseteq
\mathcal{M}\bigcap \mathcal{N}$, and the core of $W_{(\alpha ,\beta )}$ is
trivial, it is clear that all relevant weights are equal to $1$, so $%
\left\langle x,y,a\right\rangle ^{(h,\ell )}e_{(hi+m,\ell j+n)}=\left\langle
1,1,1\right\rangle e_{(hi+m,\ell j+n)}$, and therefore $\left\langle
x,y,a\right\rangle ^{(h,\ell )}|_{\mathcal{H}_{(m,n)}}\cong \left\langle
1,1,1\right\rangle $.

The proof is complete.
\end{proof}

We now recall the characterization of hyponormality, $2$-hyponormality and
subnormality for $2$-variable weighted shifts in $\mathcal{S}_{1}$ found in
\cite[Proofs of Theorems 3.1 and 3.3]{CLY4}. \ Recall that $\mathcal{S}_{1}=%
\mathcal{S}\bigcap \mathfrak{H}_{1}$, $\mathcal{S}_{2}=\mathcal{S}\bigcap
\mathfrak{H}_{2}$ and $\mathcal{S}_{\infty }=\mathcal{S}\bigcap \mathfrak{H}%
_{\infty }$.

\begin{theorem}
\label{characterization}(cf. \cite{CLY4}) \ Let $\left\langle
x,y,a\right\rangle \in \mathcal{S}_{1}$. \ Then \newline
(i) \ $\left\langle x,y,a\right\rangle \in \mathfrak{H}_{2}\iff
f_{2}(x,y,a):=(1-x^{2})-y^{2}(1-a^{2})\geq 0.$\newline
(ii) \ $\left\langle x,y,a\right\rangle \in \mathfrak{H}_{\infty }\iff
f_{2}(x,y,a)\geq 0.$
\end{theorem}

\begin{corollary}
\label{subnormal}Let $\left\langle x,y,a\right\rangle \in \mathcal{S}_{1}$.
\ The following statements are equivalent.\newline
(i) $\ \left\langle x,y,a\right\rangle \in \mathfrak{H}_{2}$;\newline
(ii) $\ \left\langle x,y,a\right\rangle \in \mathfrak{H}_{\infty }$;\newline
(iii) $\ y\leq \sqrt{\frac{1-x^{2}}{1-a^{2}}}$.
\end{corollary}

\begin{corollary}
\label{Cora}Let $0<a<1$. \ Then $\left\langle a,1,a\right\rangle \in
\mathfrak{H}_{\infty }$.

\begin{proof}
We apply Corollary \ref{subnormal} with $x:=a$ and $y:=1$. \ Since condition
(iii) is satisfied, it follows that $\left\langle a,1,a\right\rangle \in
\mathfrak{H}_{\infty }$.
\end{proof}

\begin{lemma}
\label{Lemy}Let $0<y<1$. \ Then $\left\langle 1,y,1\right\rangle \in
\mathfrak{H}_{\infty }$.

\begin{proof}
Here $T_{1}\cong I\bigotimes U_{+}$ and $T_{2}\cong S_{y}\bigotimes I$, so $%
\left\langle 1,y,1\right\rangle \equiv W_{(\alpha ,\beta )}$ is clearly
subnormal.
\end{proof}
\end{lemma}
\end{corollary}

\begin{theorem}
\label{hyponormal}Let $\left\langle x,y,a\right\rangle \in \mathcal{S}_{1}$.
\ The following statements are equivalent.\newline
(i) $\ \left\langle x,y,a\right\rangle ^{(h,\ell )}\in \mathfrak{H}_{1}$ for
all $h,\ell \geq 1$.\newline
(ii) $\ \left\langle x,y,a\right\rangle ^{(h_{0},\ell _{0})}\in \mathfrak{H}%
_{1}$ for some $h_{0},\ell _{0}\geq 1$.

\begin{proof}
It is clearly sufficient to establish (ii) $\Rightarrow $ (i). \ Assume
therefore that $\left\langle x,y,a\right\rangle ^{(h_{0},\ell _{0})}\in
\mathfrak{H}_{1}$ for some $h_{0},\ell _{0}\geq 1$. \ By Proposition \ref%
{representation}, we know that%
\begin{equation*}
\left\langle x,y,a\right\rangle ^{(h_{0},\ell _{0})}\cong \left\langle
x,y,a\right\rangle \bigoplus (h-1)\cdot \left\langle 1,\frac{ay}{x}%
,1\right\rangle \bigoplus (\ell -1)\cdot \left\langle a,1,a\right\rangle
\bigoplus (h-1)(\ell -1)\cdot \left\langle 1,1,1\right\rangle .
\end{equation*}%
An application of Corollary \ref{Cora} and Lemma \ref{Lemy} shows that $%
\left\langle x,y,a\right\rangle \in \mathfrak{H}_{1}$, \ Now, let $h,\ell
\geq 1$ be arbitrary. \ A new application of Proposition \ref{representation}
(this time using $h$ and $\ell $) shows that $\left\langle
x,y,a\right\rangle ^{(h,\ell )}\in \mathfrak{H}_{1}$. \ The proof is
complete.
\end{proof}
\end{theorem}

\begin{corollary}
Let $\left\langle x,y,a\right\rangle \in \mathcal{S}_{1}$, and let $k\geq 2$
be given. \ The following statements are equivalent.\newline
(i) $\ $For some $h_{0},\ell _{0}\geq 1$, $\left\langle x,y,a\right\rangle
^{(h_{0},\ell _{0})}\in \mathfrak{H}_{k}$.\newline
(ii) $\ $For all $h,\ell \geq 1$, $\left\langle x,y,a\right\rangle ^{(h,\ell
)}\in \mathfrak{H}_{k}$.\newline
(iii) $\ $For some $h_{0},\ell _{0}\geq 1$ $\left\langle x,y,a\right\rangle
^{(h_{0},\ell _{0})}\in \mathfrak{H}_{\infty }$.\newline
(iv) $\ $For all $h,\ell \geq 1$ $\left\langle x,y,a\right\rangle ^{(h,\ell
)}\in \mathfrak{H}_{\infty }$.
\end{corollary}

\begin{proof}
Straightforward from Proposition \ref{representation} and the Proof of
Theorem \ref{hyponormal}.
\end{proof}

We conclude this section with a problem of independent interest. \ Recall
that $\mathcal{A}=\{W_{(\alpha ,\beta )}\in \mathcal{TC}:$ the Berger
measure of $c(W_{(\alpha ,\beta )})$ is $1$-atomic$\}$.

\begin{problem}
\label{question3}Is $\mathcal{S}_{1}$ the largest class in $\mathcal{A}$ \
for which the implication
\begin{equation*}
W_{(\alpha ,\beta )}^{(h_{0},\ell _{0})}\in \mathfrak{H}_{2}\text{ for some }%
h_{0},\ell _{0}\geq 1\Rightarrow W_{(\alpha ,\beta )}\in \mathfrak{H}%
_{\infty }
\end{equation*}%
holds?
\end{problem}

%%%%%%%%%%%%%%%%%%%%%%%%

\section{Appendix}

For the reader's convenience, in this section we gather several well known
auxiliary results which are needed for the proofs of the main results in
this article. \ First, to detect hyponormality for $2$-variable weighted
shifts we use a simple criterion involving a base point $\mathbf{k}$ in $%
\mathbb{Z}_{+}^{2}$ and its five neighboring points in $\mathbf{k}+\mathbb{Z}%
_{+}^{2}$ at path distance at most $2$.

\begin{lemma}
(\cite[Theorem 6.1]{bridge})\label{joint hypo} (Six-point Test) Let $%
W_{(\alpha ,\beta )}\mathbf{\equiv (}T_{1},T_{2})$ be a $2$-variable
weighted shift, with weight sequences $\alpha $ and $\beta $. \ Then%
\begin{equation*}
\begin{tabular}{l}
$\lbrack W_{(\alpha ,\beta )}^{\ast },W_{(\alpha ,\beta )}\mathbf{]}\geq 0$
\\
$\iff H(k_{1},k_{2})(1):=\left(
\begin{array}{cc}
\alpha _{\mathbf{k}+\mathbf{\varepsilon }_{1}}^{2}-\alpha _{\mathbf{k}}^{2}
& \alpha _{\mathbf{k}+\mathbf{\varepsilon }_{2}}\beta _{\mathbf{k}+\mathbf{%
\varepsilon }_{1}}-\alpha _{\mathbf{k}}\beta _{\mathbf{k}} \\
\alpha _{\mathbf{k}+\mathbf{\varepsilon }_{2}}\beta _{\mathbf{k}+\mathbf{%
\varepsilon }_{1}}-\alpha _{\mathbf{k}}\beta _{\mathbf{k}} & \beta _{\mathbf{%
k}+\mathbf{\varepsilon }_{2}}^{2}-\beta _{\mathbf{k}}^{2}%
\end{array}%
\right) \geq 0\text{ }(\text{for all }\mathbf{k}\in \mathbb{Z}_{+}^{2})\text{%
.}$%
\end{tabular}%
\end{equation*}
\end{lemma}

Next, we present an analogous criterion for the $k$-hyponormality of $%
2$-variable weighted shifts.

\begin{lemma}
(\cite[Theorem 2.4]{CLY1})\label{khypo} \ Let $W_{(\alpha ,\beta )}\equiv
(T_{1},T_{2})$ be a $2$-variable weighted shift with weight sequence $\alpha
$ and $\beta $. \ The following statements are equivalent:\newline
(i) $\ W_{(\alpha ,\beta )}$ is $k$-hyponormal;\newline
(ii) $\ M_{\mathbf{k}}(k):=(\gamma _{\mathbf{k}+(n,m)+(p,q)})_{_{0\leq
p+q\leq k}^{0\leq n+m\leq k}}\geq 0$ for all $\mathbf{k}\in \mathbb{Z}%
_{+}^{2}$.
\end{lemma}

In particular, a commuting pair $(T_{1},T_{2})$ is $2$-hyponormal if and
only if the $5$-tuple $(T_{1},T_{2},T_{1}^{2},T_{1}T_{2},$ $T_{2}^{2})$ is
hyponormal. \ For $2$-variable weighted shifts, this is equivalent to the
condition (Fifteen-point Test)%
\begin{equation*}
M_{\mathbf{k}}(2):=(\gamma _{\mathbf{k}+(n,m)+(p,q)})_{_{0\leq p+q\leq
2}^{0\leq n+m\leq 2}}\geq 0\quad (\text{{all} }\mathbf{k}\in \mathbb{Z}%
_{+}^{2});
\end{equation*}%
that is,
\begin{equation*}
M_{\mathbf{k}}(2)\equiv \left(
\begin{array}{cccccc}
\gamma _{k_{1},k_{2}} & \gamma _{k_{1}+1,k_{2}} & \gamma _{k_{1},k_{2}+1} &
\gamma _{k_{1}+2,k_{2}} & \gamma _{k_{1}+1,k_{2}+1} & \gamma _{k_{1},k_{2}+2}
\\
\gamma _{k_{1}+1,k_{2}} & \gamma _{k_{1}+2,k_{2}} & \gamma _{k_{1}+1,k_{2}+1}
& \gamma _{k_{1}+3,k_{2}} & \gamma _{k_{1}+2,k_{2}+1} & \gamma
_{k_{1}+1,k_{2}+2} \\
\gamma _{k_{1},k_{2}+1} & \gamma _{k_{1}+1,k_{2}+1} & \gamma _{k_{1},k_{2}+2}
& \gamma _{k_{1}+2,k_{2}+1} & \gamma _{k_{1}+1,k_{2}+2} & \gamma
_{k_{1},k_{2}+3} \\
\gamma _{k_{1}+2,k_{2}} & \gamma _{k_{1}+3,k_{2}} & \gamma _{k_{1}+2,k_{2}+1}
& \gamma _{k_{1}+4,k_{2}} & \gamma _{k_{1}+3,k_{2}+1} & \gamma
_{k_{1}+2,k_{2}+2} \\
\gamma _{k_{1}+1,k_{2}+1} & \gamma _{k_{1}+2,k_{2}+1} & \gamma
_{k_{1}+1,k_{2}+2} & \gamma _{k_{1}+3,k_{2}+1} & \gamma _{k_{1}+2,k_{2}+2} &
\gamma _{k_{1}+1,k_{2}+3} \\
\gamma _{k_{1},k_{2}+2} & \gamma _{k_{1}+1,k_{2}+2} & \gamma _{k_{1},k_{2}+3}
& \gamma _{k_{1}+2,k_{2}+2} & \gamma _{k_{1}+1,k_{2}+3} & \gamma
_{k_{1},k_{2}+4}%
\end{array}%
\right) \geq 0.
\end{equation*}%
This takes into account a base point $k$ and its $14$ neighbors at path
distance at most $4$.

To check subnormality of $2$-variable weighted shifts, we introduce some
definitions. \ \newline
(i) \ Let $\mu $ and $\nu $ be two positive measures on $\mathbb{R}_{+}.$ \
We say that $\mu \leq \nu $ on $X:=\mathbb{R}_{+},$ if $\mu (E)\leq \nu (E)$
for all Borel subset $E\subseteq \mathbb{R}_{+}$; equivalently, $\mu \leq
\nu $ if and only if $\int fd\mu \leq \int fd\nu $ for all $f\in C(X)$ such
that $f\geq 0$ on $\mathbb{R}_{+}$.\newline
(ii)\ \ Let $\mu $ be a probability measure on $X\times Y$, and assume that $%
\frac{1}{t}\in L^{1}(\mu ).$ \ The \textit{extremal measure} $\mu _{ext}$
(which is also a probability measure) on $X\times Y$ is given by $d\mu
_{ext}(s,t):=(1-\delta _{0}(t))\frac{1}{t\left\Vert \frac{1}{t}\right\Vert
_{L^{1}(\mu )}}d\mu (s,t)$. \newline
(iii) \ Given a measure $\mu $ on $X\times Y$, the \textit{marginal measure}
$\mu ^{X}$ is given by $\mu ^{X}:=\mu \circ \pi _{X}^{-1}$, where $\pi
_{X}:X\times Y\rightarrow X$ is the canonical projection onto $X$. \ Thus $%
\mu ^{X}(E)=\mu (E\times Y)$, for every $E\subseteq X$. \ \newline
Then we have:

\begin{lemma}
\label{backext}\cite[Proposition 3.10]{CuYo1} \ (Subnormal Backward
Extension) \ Let $W_{(\alpha ,\beta )}$ be a $2$-variable weighted shift,
and assume that $W_{(\alpha ,\beta )}|_{\mathcal{M}_{1}}$ is subnormal with
associated measure $\mu _{\mathcal{M}_{1}}$ and that $W_{0}:=\operatorname{shift}%
\;(\alpha _{00},\alpha _{10},\cdots )$ is subnormal with associated measure $%
\xi _{0}$. \ Then $W_{(\alpha ,\beta )}$ is subnormal if and only if\newline
$(i)$ $\ \frac{1}{t}\in L^{1}(\mu _{\mathcal{M}_{1}})$;\newline
$(ii)$ $\ \beta _{00}^{2}\leq (\left\| \frac{1}{t}\right\| _{L^{1}(\mu _{%
\mathcal{M}_{1}})})^{-1}$;\newline
$(iii)$ $\ \beta _{00}^{2}\left\| \frac{1}{t}\right\| _{L^{1}(\mu _{\mathcal{%
M}_{1}})}(\mu _{\mathcal{M}_{1}})_{ext}^{X}\leq \xi _{0}$.\newline
Moreover, if $\beta _{00}^{2}\left\| \frac{1}{t}\right\| _{L^{1}(\mu _{%
\mathcal{M}_{1}})}=1,$ then $(\mu _{\mathcal{M}_{1}})_{ext}^{X}=\xi _{0}$. \
In the case when $W_{(\alpha ,\beta )}$ is subnormal, the Berger measure $%
\mu $ of $W_{(\alpha ,\beta )}$ is given by
\begin{equation*}
d\mu (s,t)=\beta _{00}^{2}\left\| \frac{1}{t}\right\| _{L^{1}(\mu _{\mathcal{%
M}_{1}})}d(\mu _{\mathcal{M}_{1}})_{ext}(s,t)+(d\xi _{0}(s)-\beta
_{00}^{2}\left\| \frac{1}{t}\right\| _{L^{1}(\mu _{\mathcal{M}_{1}})}d(\mu _{%
\mathcal{M}_{1}})_{ext}^{X}(s))d\delta _{0}(t).
\end{equation*}
\end{lemma}

\begin{lemma}
(\cite[Theorem 2.8]{Yo2})\label{thmbackext} \ Let $W_{(\alpha ,\beta )}\in
\mathfrak{H}_{0}$ be a $2$-variable weighted shift whose weight diagram is
given in Figure \ref{S1}(ii), so that $W_{(\alpha ,\beta )}\mathbf{|}_{%
\mathcal{M}_{1}}\cong (I\otimes \operatorname{shift}\;(\beta _{1},\beta _{2},\cdots
),U_{+}\otimes bI)$.\ \ Assume that $\left\| W_{\alpha }\right\| =b>0$,
where $W_{\alpha }\equiv \operatorname{shift}\;(\alpha _{0},\alpha _{1},\alpha
_{2},\cdots )$. \ Then $W_{(\alpha ,\beta )}\in \mathfrak{H}_{1}\iff
W_{(\alpha ,\beta )}\in \mathfrak{H}_{\infty }\iff $ the Berger measure $\mu
_{\alpha }$ of $W_{\alpha }$ has an atom at $b^{2}$.
\end{lemma}

Given a subnormal $2$-variable weighted shift $W_{(\alpha ,\beta )}$ with
Berger measure $\mu $, we let $W_{\alpha ^{(j)}}$$(j\geq 0)$ (resp. $%
W_{\beta ^{(i)}}\;(i\geq 0)$) denote the associated $j$-th horizontal (resp.
$i$-th vertical) slice of $W_{(\alpha ,\beta )}$. \ Clearly, $W_{\alpha
^{(j)}}$ (resp. $W_{\beta ^{(i)}}$) is subnormal, and we let $\xi _{j}$
(resp. $\eta _{i}$) denote its Berger measure. \ We proved in \cite{CuYo2}
that $d\xi _{j}(s):=\{\frac{1}{\gamma _{0j}}\int t^{j}\;d\Phi
_{s}(t)\}\;d\xi (s)$, where $d\mu (s,t)\equiv d\Phi _{s}(t)d\xi (s)$ is the
canonical disintegration of $\mu $ by vertical slices (resp. $d\eta
_{i}(t)=\{\frac{1}{\gamma _{i0}}\int s^{i}\;d\Psi _{t}(s)\}\;d\eta (t)$,
where $d\mu (s,t)\equiv d\Psi _{t}(s)d\eta (t)$ is the canonical
disintegration of $\mu $ by horizontal slices).

\begin{lemma}
(\cite[Theorem 3.3]{CuYo2})\label{necessary} \ Let $\mu $, $\xi _{j}$ and $%
\eta _{i}$ be as above. \ If $W_{(\alpha ,\beta )}\in \mathfrak{H}_{\infty }$%
, then for every $i,j\geq 0$ we have
\begin{equation}
\xi _{j+1}\ll \xi _{j}\text{ and }\eta _{i+1}\ll \eta _{i}.
\label{conditions}
\end{equation}
\end{lemma}

\begin{lemma}
\label{smu}(cf. \cite{Smu}, \cite[Proposition 2.2]{tcmp1}) \ Let $M\equiv
\left(
\begin{array}{cc}
A & B \\
B^{\ast } & C%
\end{array}%
\right) $ be a $2\times 2$ operator matrix, where $A$ and $C$ are square
matrices and $B$ is a rectangular matrix. \ Then
\begin{equation*}
M\geq 0\iff \text{there exists }W\text{ such that }\left\{
\begin{tabular}{l}
$A\geq 0$ \\
$B=AW$ \\
$C\geq W^{\ast }AW.$%
\end{tabular}%
\right.
\end{equation*}
\end{lemma}

\end{document}